\documentclass[a4paper,11pt]{amsart}
\usepackage{amssymb, graphicx, diagrams}
\newtheorem{lemma}{Lemma}
\newtheorem{prop}{Proposition}
\newtheorem{thm}{Theorem}
\newtheorem{cor}{Corollary}
\newtheorem{conj}{Conjecture}
\theoremstyle{definition}
\newtheorem{rem}{Remark}
\newtheorem{defn}{Definition}
\newtheorem{ex}{Example}

\newcounter{numl}

\newcommand{\labelnuml}{\textup{(\roman{numl})}}
\newenvironment{numlist}{\begin{list}{\labelnuml}%
{\usecounter{numl}\setlength{\leftmargin}{0pt}%
\setlength{\itemindent}{2\parindent}%
\setlength{\itemsep}{\smallskipamount}\def
\makelabel ##1{\hss \llap {\upshape ##1}}}}{\end{list}}
\newenvironment{bulletlist}{\begin{list}{\labelitemi}%
{\setlength{\leftmargin}{0pt}\setlength{\itemindent}{\parindent}%
\setlength{\itemsep}{\smallskipamount}\def
\makelabel ##1{\hss \llap {\upshape ##1}}}}{\end{list}}

\newenvironment{fact}{\smallbreak\noindent\it\ignorespaces}{\smallbreak}
\newcommand{\R}{{\mathbb R}}
\newcommand{\C}{{\mathbb C}}
\newcommand{\Z}{{\mathbb Z}}
\newcommand{\N}{{\mathbb N}}
\newcommand{\Q}{{\mathbb Q}}
\newcommand{\T}{{\mathbb T}}
\newcommand{\cP}{{\mathcal P}}
\newcommand{\cA}{{\mathcal A}}
\newcommand{\cE}{{\mathcal E}}
\newcommand{\cI}{{\mathcal I}}
\newcommand{\cL}{{\mathcal L}}
\newcommand{\cK}{{\mathcal K}}
\newcommand{\cO}{{\mathcal O}}
\newcommand{\vE}{E}
\newcommand{\lE}{{\mathcal E}}
\newcommand{\zE}{e}
\newcommand{\mult}{^{\scriptscriptstyle\times}}

\newcommand{\Scal}{\mathit{Scal}}
\newcommand{\Id}{\mathit{Id}}
\newcommand{\eps}{\varepsilon}

\newcommand{\trace}{\mathop{\mathrm{tr}}\nolimits}
\newcommand{\pfaff}{\mathop{\mathrm{pf}}\nolimits}
\newcommand{\grad}{\mathop{\mathrm{grad}}\nolimits}
\newcommand{\vspan}{\mathop{\mathrm{span}}\nolimits}
\newcommand{\rk}{\mathop{\mathrm{rank}}\nolimits}

\newcommand{\Hom}{\mathrm{Hom}}

\newcommand{\Aut}{\mathrm{Aut}}
\newcommand{\Vol}{\mathrm{Vol}}
\newcommand{\pr}{\mathrm{pr}}
\newcommand{\restr}[1]{|_{#1}^{\vphantom x}}
\newcommand{\Mp}{p}
\newcommand{\Mpc}{p_{\mathrm{c}}}
\newcommand{\Mpn}{p_{\mathrm{nc}}}
\newcommand{\cb}{\overline{c}}
\newcommand{\mapsinto}{\hookrightarrow}
\begin{document}
\title[Hamiltonian 2-forms in K{\smash{\"a}}hler geometry, III]
{Hamiltonian 2-forms in K{\smash{\"a}}hler geometry,\\
III Extremal metrics and stability}
\author[V. Apostolov]{Vestislav Apostolov}
\address{Vestislav Apostolov \\ D{\'e}partement de Math{\'e}matiques\\
UQAM\\ C.P. 8888 \\ Succ. Centre-ville \\ Montr{\'e}al (Qu{\'e}bec) \\
H3C 3P8 \\ Canada}
\email{apostolo@math.uqam.ca}
\author[D. Calderbank]{David M.~J.~Calderbank}
\address{David M. J. Calderbank \\ Department of Mathematics\\
University of York\\ Heslington\\ York YO10 5DD\\ England}
\email{dc511@york.ac.uk}
\author[P. Gauduchon]{Paul Gauduchon}
\address{Paul Gauduchon \\ Centre de Math\'ematiques\\
{\'E}cole Polytechnique \\ UMR 7640 du CNRS
\\ 91128 Palaiseau \\ France}
\email{pg@math.polytechnique.fr}
\author[C. T\o nnesen-Friedman]{Christina W.~T\o nnesen-Friedman}
\address{Christina W. T\o nnesen-Friedman\\
Department of Mathematics\\ Union College\\
Schenectady\\  New York  12308\\ USA }
\email{tonnesec@union.edu}
\thanks{The first author was supported by NSERC grant OGP0023879, the second
author by an EPSRC Advanced Research Fellowship and the fourth author by the
Union College Faculty Research Fund.}
\date{\today}
\begin{abstract}  This paper concerns the explicit construction of extremal
K\"ahler metrics on total spaces of projective bundles, which have been
studied in many places. We present a unified approach, motivated by the theory
of hamiltonian $2$-forms (as introduced and studied in previous papers in the
series) but this paper is largely independent of that theory.

We obtain a characterization, on a large family of projective bundles, of
those `admissible' K\"ahler classes (i.e., the ones compatible with the bundle
structure in a way we make precise) which contain an extremal K\"ahler
metric. In many cases, such as on geometrically ruled surfaces, every K\"ahler
class is admissible. In particular, our results complete the classification of
extremal K\"ahler metrics on geometrically ruled surfaces, answering several
long-standing questions.

We also find that our characterization agrees with a notion of K-stability
for admissible K\"ahler classes. Our examples and nonexistence results
therefore provide a fertile testing ground for the rapidly developing theory
of stability for projective varieties, and we discuss some of the
ramifications. In particular we obtain examples of projective varieties which
are destabilized by a non-algebraic degeneration.
\end{abstract}
\maketitle

In this paper we give a systematic overview of some explicit constructions of
extremal K\"ahler metrics on projective bundles and relate our constructions
to the theory of stability for algebraic varieties. Recall that a K\"ahler
metric $(g,J,\omega)$ on a manifold $M$ is said to be
{\it extremal}~\cite{calabi1} if the scalar curvature is a {\it Killing
potential}, i.e., its symplectic gradient is a Killing vector field, or
equivalently, its gradient is a holomorphic vector field. If $(M,J)$ is
compact, the extremal K\"ahler metrics are exactly the critical points of the
$L_2$-norm of the scalar curvature on the space of compatible K\"ahler metrics
in a fixed K\"ahler class $\Omega=[\omega]$.

There has been a great deal of interest recently in the relation between
extremal K\"ahler metrics, especially constant scalar curvature (CSC) K\"ahler
metrics, and {\it stability}: roughly speaking there are conjectures that the
existence of an extremal or CSC K\"ahler metric in an integral K\"ahler class
$\Omega$ on a compact complex manifold $M$ should be equivalent to an
algebraic geometric stability condition on the Kodaira embeddings of $M$ into
$P(H^0(M,L^k)^*)$ for $k\gg 1$ where $c_1(L)=\Omega/2\pi$.

Our own interest in constructions of extremal K\"ahler metrics has been
renewed and stimulated not only by these developments, but also by a unifying
principle, which we shall explain and apply here, underlying explicit examples
of such metrics on projective bundles obtained
in~\cite{calabi2,guan0,hwang,hwang-singer,christina1}.  In short these
examples have in common the presence of a {\it hamiltonian $2$-form of order
$1$}. 

A {\it hamiltonian $2$-form}~\cite{ACG2} on a K\"ahler manifold
$(M,J,g,\omega)$, of real dimension $2m>2$ is a real $(1,1)$-form (i.e., a
$J$-invariant $2$-form) $\phi$ such that
\begin{equation*}
2 \nabla _{X} \phi = d \trace\phi \wedge (JX)^{\flat} - (J d \trace\phi ) 
\wedge X^{\flat}
\end{equation*}
for all $X \in TM$ (where $X^{\flat}(Y) = g(X,Y)$ for $Y \in TM$ and
$\trace\phi = \langle\omega, \phi\rangle_{g}$).  The {\it momentum polynomial}
of $\phi$ is then defined to be
\begin{equation*}
p(t) :=  (-1)^m\pfaff (\phi-t\omega)
=  t^m - (\trace\phi) \,t^{m-1} + \cdots + (-1)^m \pfaff\phi,
\end{equation*}
where the {\it pfaffian} is defined by $\phi \wedge \cdots \wedge \phi =
(\pfaff\phi)\omega\wedge\cdots\wedge\omega$. The reason $\phi$ is called
{\it hamiltonian} is that for any $t\in\R$, $\Mp(t)$ is a hamiltonian for
a Killing vector field $K(t) := J \grad_g \Mp(t)$ (this is not difficult to
show~\cite[\S2]{ACG2}).  The integer $\ell=\max_{x\in M} \dim\vspan
\{K(t)_x:t\in\R\}$ (with $0\leq \ell\leq m$) is called the {\it order} of
$\phi$.

We do not wish to impose the study of hamiltonian $2$-forms on the reader of
this paper, since we only need the general theory as motivation for the
classes of complex manifolds and metrics that we shall study. We therefore now
state a classification result which reduces the theory of hamiltonian
$2$-forms of order $1$ to an Ansatz for metrics on projective bundles.  This
result follows easily from~\cite{ACG2, ACGT3}, as we explain in Appendix A,
where we also explain why we restrict attention to the order $1$ case.

\begin{thm}\label{ACGTthm}
Let $(M,g,J,\omega)$ be a compact connected K\"ahler $2m$-manifold with a
hamiltonian $2$-form $\phi$ of order $1$. Then\textup:
\begin{bulletlist}
\item there is an effective isometric hamiltonian $S^1$ action on $M$
generated by a vector field $K=J\grad_g z$ such that the stable quotient of
$M$ by the induced holomorphic $\C\mult$ action is a compact connected complex
manifold $\smash{\hat S}$ of real dimension $2(m-1)$\textup;

\item without loss, the image of the momentum map $z$ is $[-1,1]$, and there
are K\"ahler manifolds $S_a$ of dimension $2d_a$ and real numbers $x_a$,
indexed by $a\in\smash{\hat\cA}\subset\N\cup\{\infty\}$, such that $\hat{S}$
is covered by $\prod_{a\in\hat\cA} S_a$ and $0 < |x_{a}| \leq 1$ with equality
iff $a\in\{0,\infty\}$\textup;

\item $z$ is a Morse--Bott function~\cite{atiyah} on $M$ with critical set
$z^{-1}(\{-1,1\})$, $M^0:=z^{-1}((-1,1))$ is a principal $\C\mult$-bundle over
$\smash{\hat S}$, and on each $S_a$ there is a K\"ahler metric $(\pm g_a,\pm
\omega_a)$, with the sign chosen so that $\omega_a/x_a$ is positive, such that
the K\"ahler structure on $M^0$ is
\begin{equation}\label{metric}
g=\sum_{a\in\smash{\hat\cA}\vphantom{I}} \frac{1+x_a z}{x_a} g_a
+\frac{dz^{\smash 2}}{\Theta(z)}+\Theta(z)\theta^2,\quad
\omega=\sum_{a\in\smash{\hat\cA}\vphantom{I}} \frac{1+x_a z}{x_a}\omega_a
+dz\wedge \theta,
\end{equation}
where $\theta$ is a connection $1$-form \textup($\theta(K)=1$\textup) with
$d\theta=\sum_{a\in\smash{\hat\cA}\vphantom{I}} \omega_a$, and $\Theta$ is a
smooth function on $[-1,1]$ satisfying
\begin{gather}\label{positivity}
\Theta>0\quad\text{on}\quad (-1,1),\\
\Theta(\pm 1) = 0,\qquad
\Theta'(\pm 1) = \mp 2;
\label{boundary}
\end{gather}
\item if $0\in\smash{\hat\cA}$ then $x_0=1$, $S_0 = \C P^{d_0}$ and $g_0$ is
the Fubini--Study metric of scalar curvature $2 d_0 (d_0 + 1)$,
otherwise we set $d_0=0$\textup; likewise, if $\infty\in\smash{\hat\cA}$ then
$x_\infty=-1$, $S_\infty = \C P^{d_\infty}$ and $-g_\infty$ is the
Fubini--Study metric of scalar curvature $2 d_\infty (d_\infty + 1)$,
otherwise we set $d_\infty=0$\textup; we also put
$\cA=\smash{\hat\cA}\smallsetminus \{0,\infty\} \subset\Z^+$\textup;

\item the blow-up $\smash{\hat M}$ of $M$ along $z^{-1}(\{-1,1\})$ is
$\C\mult$-equivariantly biholomorphic to $M^0\times_{\C\mult} \C P^1 \to
\hat{S}$ and $\hat{S}$ is a fibre product of flat projective unitary $\C
P^{d_0}$- and $\C P^{d_\infty}$-bundles over a K\"ahler manifold $S$ covered
by $\prod_{a\in\cA} S_a$.
\end{bulletlist}
\underline{If we assume that} $\smash{\hat S} = P(\vE_0)\times_S
P(\vE_\infty)\to S$ for projectively-flat hermitian vector bundles
$\vE_0,\vE_\infty\to S$, then these bundles can be chosen so that $M$ is
$\C\mult$-equivariantly biholomorphic to $P(\vE_0\oplus \vE_\infty)\to S$, and
we therefore have $\cb_1(\vE_\infty)-\cb_1(\vE_0) = \sum_{a\in \cA}
[\omega_a/2\pi]$, where $\cb_1(E)= c_1(E)/\rk E$.
\end{thm}
The final assumption of this theorem is automatic if $\pi_1(S)=1$, when
$\smash{\hat S}= \C P^{d_0}\times S\times \C P^{d_\infty}$ and there is a line
bundle $\cL\to S$ with $c_1(\cL)=\sum_{a\in \cA}[\omega_a/2\pi]$ such that $M$
is $\C\mult$-equivariantly biholomorphic to $P(\cO\otimes\C^{d_0+1}
\oplus\cL\otimes\C^{d_\infty+1} )\to S$.

We shall use this theorem as an Ansatz for constructing extremal K\"ahler
metrics on projective bundles of the form $P(\vE_0\oplus \vE_\infty)\to S$.
For this we shall use the following elementary
computation~\cite{guan0,hwang,ACG2}.
\begin{prop}\label{equations}
Let $g$ be a K\"ahler metric of the form~\eqref{metric} on $M^0$ and write
$F(z)=\Theta(z)\Mpc(z)$ with $\Mpc(z):=\prod_{a\in\hat\cA} (1+x_az)^{d_a}$.
Then $g$ is extremal, with $\Scal_g$ a constant affine function of $z$, iff
\begin{bulletlist}
\item there is a polynomial $P$ of degree $\leq \#\smash{\hat\cA}+1$ such that
\begin{equation}\label{extremal-deriv}
F''(z) = \biggl( \prod_{a\in\smash{\hat\cA}\vphantom{I}} (1+x_a z)^{d_a-1}
\biggr) P(z);
\end{equation}
\item for all $a$, $\pm g_a$ has constant scalar curvature $\Scal_{\pm
g_a}=\pm 2d_as_a$ where
\begin{equation}\label{extremal-values}
P(-1/x_a) = 2d_as_a x_a \prod_{b\in\smash{\hat\cA}\setminus\{a\}}
\biggl(1-\frac{x_b}{x_a}\biggr).
\end{equation}
\end{bulletlist}
The metric $g$ then has constant scalar curvature iff $P$ has degree $\leq
\#\smash{\hat\cA}$.
\end{prop}
Compared to~\cite{ACG2,ACGT3}, we have rescaled $F(z)$ and $\Mpc(z)$ by
$\prod_{a\in\hat\cA} x_a^{d_a}$: this is convenient as $\Mpc(z)$ is then
positive on $(-1,1)$. Thus $\Theta$ is positive on $(-1,1)$ if and only if $F$
is. Also if $\Theta(z)$ satisfies~\eqref{boundary}, then
$F(z)=\Theta(z)\Mpc(z)$ satisfies
\begin{equation}\label{BCs}
F(\pm1)=0,\qquad F'(\pm 1) = \mp 2 \Mpc(\pm 1).
\end{equation}

The structure of the paper is as follows. In section~\ref{admissible} we study
metrics of the form~\eqref{metric} and show that the conditions of
Theorem~\ref{ACGTthm} are sufficient for the compactification of such metrics
on a projective bundle $M=P(\vE_0\oplus\vE_\infty)\to S$.  We shall call
metrics of this form (compatible with the given projective bundle structure on
$M$ and local product structure on $S$) up to scale {\it admissible}.  The
projective bundle $M$ is similarly called an {\it admissible bundle} or
{\it manifold}.  We describe the {\it admissible K\"ahler classes} (i.e.,
those containing an admissible metric) in~\S\ref{admissible-class} and
parameterize the admissible metrics in a given admissible class
in~\S\ref{compactify}. We end the section by computing in~\S\ref{isometry},
the Lie algebra of Killing vector fields with zeros on $M$.

In section~\ref{main} we study extremal K\"ahler metrics in admissible
K\"ahler classes and prove our first main result. In preparation for this, we
consider what implications the existence of an extremal K\"ahler metric has
for admissible bundles, K\"ahler classes and metrics. In
\S\ref{matsushima-lichnerowicz}, we apply the Matsushima--Lichnerowicz
criterion~\cite{lichne,matsushima} to obtain information about the
automorphism group of an admissible bundle when the base metric on $S$ is
CSC. Then, in \S\ref{futaki-invariant} we compute the Futaki invariant
$\mathfrak F_\Omega(K)$ of the vector field $K=J\grad_g z$ for any admissible
K\"ahler class $\Omega$ on an admissible $2m$-manifold $M$.  Finally in
\S\ref{K-energy-ext} we give a formula for the Mabuchi--Guan--Simanca
(modified) K-energy functional on admissible K\"ahler metrics in $\Omega$: it
is determined by a polynomial $F_\Omega$ of degree $\leq m+2$, which we call
the {\it extremal polynomial} (the degree $m+2$ coefficient being a nonzero
multiple of $\mathfrak F_\Omega(K)$).  Using the recent work of
Chen--Tian~\cite{CT,CT2}, we deduce that if $\Omega$ contains an extremal
K\"ahler metric, then $F_\Omega$ must be nonnegative on $[-1,1]$.

In \S\ref{char}, we apply Proposition~\ref{equations} to give a construction
of admissible extremal K\"ahler metrics, unifying and generalizing work of
Calabi, Guan, Hwang, Hwang--Singer and the fourth
author~\cite{calabi2,guan0,hwang,hwang-singer,christina1}. This leads to the
following result.
\begin{thm}\label{main-thm} Let $M=P(\vE_0\oplus\vE_\infty)\to S$ be an
admissible manifold, where the base $S$ is a local K\"ahler product of CSC
metrics $(\pm g_a, \pm\omega_a)$. Then there is an extremal K\"ahler metric in
an admissible K\"ahler class $\Omega$ if and only if the extremal polynomial
$F_\Omega$ is positive on $(-1,1)$. This condition always holds if $\Omega$ is
`sufficiently small'\textup; if it does, there is an admissible extremal
K\"ahler metric in $\Omega$, which is CSC if and only if the Futaki invariant
$\mathfrak F_\Omega(K)$ vanishes \textup(i.e., $F_\Omega$ has degree $\leq
m+1$\textup).

The admissible K\"ahler classes containing an extremal K\"ahler metric form a
nonempty open subset of all such classes, and those containing a CSC K\"ahler
metric form a real analytic hypersurface which is nonempty if
$\cb_1(\vE_\infty)-\cb_1(\vE_0)$ is strictly indefinite \textup(i.e., the
K\"ahler forms $\pm \omega_a$ do not all have the same sign\textup).
\end{thm}
Here we say that an admissible K\"ahler class is {\it sufficiently small} if
the $x_a$ ($a\in\cA$) are sufficiently small. Geometrically, this means that
the base $S$ is large (low curvature) compared to the fibres (high curvature).
Thus the above theorem asserts the existence of extremal K\"ahler metrics with
curvature concentrated in the fibres (cf.~\cite{hong,RT}). Note also that
`the' extremal K\"ahler metric in $\Omega$ (if it exists) is unique up to
automorphism by Chen--Tian~\cite{CT,CT2}.

In section~\ref{existence?} we present further existence and nonexistence
results for extremal and CSC metrics by computing the extremal polynomial on
various examples and testing its positivity on $(-1,1)$. In many of these
examples, every K\"ahler class on $M$ is admissible (see
Remark~\ref{classes-adm}), and therefore Theorem~\ref{main-thm} describes
exactly which K\"ahler classes contain an extremal K\"ahler metric. This is
the case when $M=P(\cO\oplus\cL)\to S$ is an admissible geometrically ruled
complex surface; the extremal polynomial $F_\Omega(z)$ is then a quartic
divisible by $1-z^2$. We thus obtain a complete resolution of the existence
question for extremal K\"ahler metrics on these complex surfaces: the K\"ahler
cone is a cone on an open interval $(a,b)$; the extremal K\"ahler metrics are
precisely those of~\cite{calabi1,hwang,christina1}, which are admissible and
locally cohomogeneity one, with K\"ahler classes parameterized by a cone on a
subinterval $(a,c)$ and $c=b$ if and only if $S$ has genus $0$ or $1$.  As
observed in~\cite{AT}, this fills in the missing step in the complete
classification of extremal K\"ahler metrics on geometrically ruled complex
surfaces.

Our results provide a fertile testing ground for the conjectures relating
extremal and CSC K\"ahler metrics to stability, and we explore this in
section~\ref{stab}. In~\S\ref{stab-CSC} we relate our results to those of
Ross--Thomas~\cite{RT} and Hong~\cite{hong}: in particular, we show that there
are CSC metrics on projective bundles $P(E)\to S$ for which $E$ is only
(slope) polystable with respect one K\"ahler class on $S$ up to scale.

In~\S\ref{ext-exist}, we relate Theorem~\ref{main-thm} to the notions of
K-polystability~\cite{Do2,tian} or relative K-polystability~\cite{szekelyhidi}
for K\"ahler classes, which are conjectured to be equivalent to the existence
of a CSC or extremal K\"ahler metric in a given class. Actually, to be
precise, we use a closely related notion of (relative) slope K-polystability
suggested by the work of Ross--Thomas~\cite{RT,RT2}.  Then, generalizing a
calculation of G.~Sz\'ekelyhidi for ruled surfaces, we establish the following
result.

\begin{thm}\label{stab=>ext} Let $\Omega$ be an admissible integral K\"ahler
class on $M=P(E_0\oplus E_\infty)\to S$, where $S$ is CSC. If $\Omega$ is
slope K-polystable, it contains a CSC K\"ahler metric, and if it is slope
K-polystable relative to $K=J\grad_g z$, it contains an extremal K\"ahler
metric.
\end{thm}
It is natural to ask if (relative) K-polystability in the sense
of~\cite{Do2,szekelyhidi} implies the existence of a CSC (or extremal)
K\"ahler metric in $\Omega$. We find that this is true if $\dim S\leq 4$, but
for $\dim S\geq 6$, we are only able to show that the extremal polynomial is
positive on $(-1,1)\cap\Q$.

Before our work, it was believed that K-polystability implies slope
K-polystability in general, but the proof in~\cite{RT,RT2} only shows that it
implies slope K-semistability, the gap being closely related to the issue of
positivity (versus nonnegativity) of the extremal polynomial at irrational
points in $(-1,1)$.

To show that this is a genuine problem, we end with some examples, on
projective line bundles over a product of three Riemann surfaces, of integral
admissible K\"ahler classes $\Omega$ such that $F_\Omega$ is positive on
$(-1,1)\cap\Q$ but has an irrational repeated root in $(-1,1)$.  We find these
examples intriguing, since by Theorem~\ref{main-thm}, these K\"ahler classes
do {\it not} contain an extremal K\"ahler metric so they should be unstable.
However, despite being projective varieties, the degeneration that
demonstrates this instability is not algebraic. While we cannot prove that
there is no other (algebraic) test configuration which would detect this
instability, it is difficult to imagine how such a test configuration could be
constructed. Our results then suggest that the non-algebraic degenerations
implicit in the use of slope K-polystability may be essential to relate
stability to existence of CSC and extremal K\"ahler metrics.

We would like to thank Claude LeBrun for helpful comments and Richard Thomas
for useful discussions concerning~\cite{RT,RT2}.

\tableofcontents

\section{Admissible bundles and K\"ahler metrics}\label{admissible}

\subsection{Admissible projective bundles}
\label{projbundle}

We use Theorem~\ref{ACGTthm} (including the final assumption) as motivation
for the class of compact complex manifold we will study.  A projective bundle
of the form $M = P\bigl(\vE_0 \oplus \vE_\infty \bigr) \stackrel{p}{\to} S$
will be called {\it admissible} or an {\it admissible manifold} if:
\begin{bulletlist}
\item $S$ is a covered by a product $\smash{\tilde S}=\prod_{a\in \cA}
S_a$ (for $\cA\subset\Z^+$) of simply-connected K\"ahler manifolds
$(S_a,\pm g_a,\pm \omega_a)$ of real dimensions $2d_a$;
\item $\vE_0$ and $\vE_{\infty}$ are holomorphic projectively-flat hermitian
vector bundles over $S$ of ranks $d_0+1$ and $d_\infty+1$ with
$\cb_1(\vE_\infty)-\cb_1(\vE_0) = [\omega_S/2\pi]$ and $\omega_S=\sum_{a\in
\cA}\omega_a$.
\end{bulletlist}
The second condition (cf.~\cite{kob}) means that we can fix hermitian metrics
on $E_0$ and $E_\infty$ whose Chern connections have tracelike curvatures
$\Omega_0\otimes\Id_{E_0}$ and $\Omega_\infty\otimes\Id_{E_\infty}$ satisfying
$\Omega_\infty-\Omega_0=\sum_{a\in\cA} \omega_a$. We normalize the induced
fibrewise Fubini--Study metrics $(g_0,\omega_0)$ and
$(-g_\infty,-\omega_\infty)$ on $P(E_0)$ and $P(E_\infty)$ to have scalar
curvatures $2d_0(d_0+1)$ and $2d_\infty(d_\infty+1)$.

We collect a few remarks and notations that we will use. We omit pullbacks by
obvious projections in these remarks.
\begin{numlist}
\item We sometimes let the index $a$ take values in $\N\cup\{\infty\}$ by
setting $d_a=0$ for $a\notin\cA\cup\{0,\infty\}$ (so that $S_a$ is a point and
$\omega_a=0$). This range will be assumed unless otherwise stated. We set
$\smash{\hat\cA}:= \{a:d_a>0\}$ so that $\cA=\smash{\hat\cA}\cap\Z^+$.
\item The pullbacks of $\vE_0$ and $\vE_\infty$ to $\smash{\tilde S}$ are of
the form $\lE_0\otimes\C^{d_0+1}$ and $\lE_\infty\otimes\C^{d_\infty+1}$,
where $\cL:=\lE_0^{-1}\otimes\lE_\infty=\bigotimes_{a\in\cA}\cL_a$ for line
bundles $\cL_a\to S_a$ with $c_1(\cL_a)=[\omega_a/2\pi]$.
\item $\zE_0:=P(\vE_0\oplus 0)$ and $\zE_\infty:= P(0\oplus \vE_\infty)$
denote the `zero' and `infinity' subbundles of $M$, covered by $S_0\times
\smash{\tilde S}$ and $\smash{\tilde S}\times S_\infty$, where $S_0=\C
P^{d_0}$ and $S_\infty=\C P^{d_\infty}$.
\item The blow-up of $M$ along $\zE_0\cup \zE_\infty$ is $\smash{\hat
M}:=P(\cO\oplus \smash{\hat\cL})\stackrel{\hat p}{\to} \smash{\hat S}$, where
$\smash{\hat S}= P(\vE_0)\times_S P(\vE_\infty) \to S$ and
$\smash{\hat\cL}=\cO(1)_{\vE_0}\otimes\cO(-1)_{\vE_\infty}$, using the
(fibrewise) hyperplane and tautological line bundles; we have
$c_1(\smash{\hat\cL}) = [\omega_{\hat S}/2\pi]$, where $\omega_{\hat
S}=\sum_a\omega_a$. If $d_0>0$ or $d_\infty>0$ we say {\it a blow-down
occurs}.
\item $\smash{\hat\zE}_0$ and $\smash{\hat\zE}_\infty$ denote the zero and
infinity sections of $\smash{\hat M}$. The pullback of $\smash{\hat\cL}$ to
$S_0\times \smash{\tilde S}\times S_\infty$ is
$\cL_0\otimes\cL\otimes\cL_\infty$, where $\cL_0=\cO(1)\to S_0$ and
$\cL_\infty=\cO(-1)\to S_\infty$.
\item $\smash{\hat S}$ has a family of local K\"ahler product metrics $g_{\hat
S}(z)$ with K\"ahler forms $z\omega_{\hat S}+\sum_a\omega_a/x_a$ and we set
$g_{\hat S}=g_{\hat S}(0)$. (Note that $g_{\hat S}$ is not compatible with
$\omega_{\hat S}$---the latter is symplectic, but not a K\"ahler form in
general.) We let $g_S(z)$ and $g_S=g_S(0)$ denote the induced local K\"ahler
product metrics on $S$.
\end{numlist}

We summarize the set-up with the following diagram of bundles and a blow-up:
\begin{diagram}[size=1.5em]
\smash{\hat M}= P(\cO\oplus \smash{\hat\cL}) &\rTo&
\smash{\hat S}=P(\vE_0)\times_S P(\vE_\infty)\\
\dTo & & \dTo\\
M= P(E_0\oplus E_\infty)&\rTo& S,
\end{diagram}
the universal cover (omitting pullbacks) of this diagram being:
\begin{diagram}[size=1.5em]
P(\cO\oplus \smash{\hat\cL}) &\rTo&
\C P^{d_0}\times \smash{\tilde S}\times \C P^{d_\infty}\\
\dTo & & \dTo\\
P(\cO\otimes \C^{d_0+1}\oplus \cL\otimes \C^{d_\infty+1})
&\rTo& \tilde S ={\textstyle\prod_{a\in\cA} S_a}.
\end{diagram}

\begin{rem}\label{integrality}

The existence of the line bundle $\smash{\hat\cL}\to\smash{\hat S}$ with
$c_1(\smash{\hat \cL})=[\omega_{\hat S}/2\pi]$ implies that $\omega_{\hat S}$
is {\it integral} in the sense that $[\omega_{\hat S}/2\pi]$ is in the image
of $H^2(\smash{\hat S},\Z)$ in $H^2(\smash{\hat S},\R)$. When $\smash{\hat S}$
is a global K\"ahler product (i.e., we can write
$M=P(\cO\otimes\C^{d_0+1}\oplus\cL\otimes \C^{d_\infty+1})\to
S=\prod_{a\in\cA}S_a$) this integrality condition means that each $\omega_a$
is integral, i.e., the compact manifolds $(S_a,\pm g_a,\pm\omega_a)$ are
Hodge.

We write $\omega_a=q_a \alpha_a$ for an integer $q_a\neq 0$, where $\alpha_a$
is a primitive integral K\"ahler form on $S_a$, so that $q_a$ is a nonzero
integer with the same sign as $(g_a,\omega_a)$, and $q_0=1$ and $q_\infty=-1$.
We now compare $[\omega_a/2\pi]$ to the first Chern class
$c_1(\cK_a^{-1})=[\rho_a/2\pi]$ of the anticanonical bundle of $S_a$, by
writing $[\rho_a] = p_a[\alpha_a]+[\rho_a]_0$, for a rational number $p_a$,
where $[\rho_a]_0\cdot[\alpha_a]^{d_a-1}=0$.  Since any line bundle $\cP$ with
first Chern class $[\alpha_a/2\pi]$ is ample, $\cP^{d_a+1}\otimes \cK_a$ is
nef by a result of Fujita~\cite{fujita} (see also \cite[Theorem
8.3]{demailly}), from which it follows easily that $p_a\leq d_a+1$.  If $S_a$
is a Riemann surface of genus ${\mathbf g}_a$, then $p_a=2(1-{\mathbf
g}_a)\leq 2$.

We set $s_a=p_a/q_a$. When $\pm g_a$ is CSC, we have $\Scal_{\pm g_a}=\pm 2
d_a s_a$, where the sign is that of $q_a$, so the scalar curvature of $\pm
g_a$ has the same sign as $p_a$. Thus, in the case of a CSC Hodge manifold
$S_a$, the Fujita inequality $p_a\leq d_a+1$ is (since $|q_a|\geq 1$)
equivalent to $\Scal_{\pm g_a} \le 2d_a(d_a+1)$.
\end{rem}

The conditions of Theorem~\ref{ACGTthm} are also sufficient for the
compactification of metrics of the form \eqref{metric} on an admissible
projective bundle $M=P(\vE_0\oplus\vE_\infty)\to S$, where $z\colon M\to
[-1,1]$ with $\zE_0=z^{-1}(1)$ and $\zE_\infty=z^{-1}(-1)$, and $\theta$ is a
connection $1$-form. Before discussing this, we introduce the K\"ahler classes
to which they belong.

\subsection{Admissible K\"ahler classes and canonical metrics}
\label{admissible-class}

Suppose that $M=P(E_0\oplus E_\infty)\to S$ is an admissible bundle.  We say
that a K\"ahler class $\Omega$ on $M$ is {\it admissible} if there are
constants $x_a$, with $x_0=1, x_\infty=-1$, such that the pullback of
$\Omega$ to $\smash{\hat M}$ has the form
\begin{equation*}
\sum_a [\omega_a]/x_a+\smash{\hat\Xi}
\end{equation*}
up to scale, where the $2$-forms $\omega_a$ are viewed as pullbacks to
$\smash{\hat M}$ of the corresponding forms on $\smash{\hat S}$ (induced by
the local product K\"ahler structure $\prod_a S_a$) and $\smash{\hat\Xi}$ is
Poincar\'e dual to $2\pi[\hat e_0+\hat e_\infty]$. Thus $\smash{\hat\Xi}=2\pi
c_1(V\smash{\hat M})$, where $V\smash{\hat M}=\cO(2)_{\cO\oplus\hat\cL}\otimes
\hat p^*\smash{\hat\cL}$ and $\cO(-1)_{\cO\oplus\hat\cL}$ is the (fibrewise)
tautological bundle of $\smash{\hat M}=P(\cO\oplus\smash{\hat\cL})$.  (The
first Chern class $[\omega_{\hat S}/2\pi]$ of $\smash{\hat\cL}$ itself pulls
back to $\smash{\hat M}$ to give the Poincar\'e dual of
$[\smash{\hat\zE}_0-\smash{\hat\zE}_\infty]$.)

It follows that admissible K\"ahler classes have the form
\begin{equation*}
\Omega=\sum_{a\in\cA} [\omega_a]/x_a+\Xi
\end{equation*}
up to scale, where the pullback of $\Xi$ to $\smash{\hat M}$ is
$[\omega_0]-[\omega_\infty]+\smash{\hat\Xi}$. Since pullback to a blow-up is
injective on cohomology, admissible K\"ahler classes on $M$ are uniquely
determined by the parameters $x_a$.

If $(g,\omega)$ is any K\"ahler metric on $M$ of the form~\eqref{metric} on
$M^0$, then we claim $\Omega=[\omega]$ is admissible. For this we first note
that on $M^0$, the K\"ahler form $\omega$ is a linear combination
$\sum_{a\in\cA} \omega_a /x_a + \eta$, where
\begin{equation*}
\eta=(z+1)\omega_0 + \sum_{a\in\cA} z \omega_a + (z-1)\omega_\infty
+ dz\wedge\theta.
\end{equation*}
Here $\omega_0$, $\omega_\infty$ and $\theta$ are defined only on $M^0$.
However, for $a\in \cA$, $\omega_a$ extends to a closed $2$-form on $M$ (as a
pullback from $S$), so $\eta$ is globally defined and closed on $M$ (since
$\omega$ is). The pullback of $\eta$ to $\smash{\hat M}$ may be written
$\omega_0-\omega_\infty+\hat\eta$ with $\hat\eta=d(z\theta)$ on $M^0$, and
since $\omega_0$ and $\omega_\infty$ are well-defined and closed on
$\smash{\hat M}$ (as pullbacks from $\smash{\hat S}$), so is $\hat\eta$, and
we easily see\footnote{On each fibre of $\hat p\colon\smash{\hat M}\to
\smash{\hat S}$, $\hat\eta/4\pi$ integrates to $1$ and so $[\hat\eta/4\pi]$
restricts to give the generator of $H^2(\hat p^{-1}(x),\Z)$. Hence by the
Leray--Hirsch theorem, $H^2(\smash{\hat M},\R)$ is generated by $[\hat\eta]$
and pullbacks from $S$. The restriction of $[\hat \eta/2\pi]$ to
$\smash{\hat\zE}_0$ is the first Chern class $[\omega_{\hat S}/2\pi]$ of
$\smash{\hat\cL}$ (and the restriction to $\smash{\hat\zE}_\infty$ is the
first Chern class $[-\omega_{\hat S}/2\pi]$ of $\smash{\hat\cL}^{-1}$). Thus
$[\hat\eta/4\pi]$ is a projective version of the Thom class of a vector
bundle.} that $[\hat \eta]=\smash{\hat\Xi}$.

Observe that $\eta$ depends implicitly on the choice of metric $(g,\omega)$ on
$M$ because the momentum map $z$ does. However, the above shows that the
cohomology class $[\eta]$ is $\Xi$, independent of this choice. From this
realisation of $\Xi$ it follows easily, by pulling back to $\zE_0$ and
$\zE_\infty$, that for a cohomology class of the form $\sum_{a\in\cA}
[\omega_a]/x_a+\Xi$ to be a K\"ahler class, it is necessary that for $a\in
\cA$, $0<|x_a|<1$ with the sign of $x_a$ such that $\omega_a/x_a$ is positive.
Conversely, we claim that any cohomology class of this form (with $0<|x_a|<1$
and $\omega_a/x_a$ positive for $a\in\cA$) is an admissible K\"ahler class and
contains a K\"ahler metric of the form~\eqref{metric} on $M^0$ up to scale. To
do this we construct a distinguished K\"ahler metric in each such class.

Let $r_0$ and $r_\infty$ be the norm functions induced by the hermitian
metrics on $E_0$ and $E_\infty$.  Then $z_0=\frac12 r_0^2$ and
$z_\infty=\frac12 r_\infty^2$ are fibrewise momentum maps for the $U(1)$
actions given by scalar multiplication in $E_0$ and $E_\infty$, generated by
$K_0$ and $K_\infty$. We equip $M$ with a fibrewise Fubini--Study metric
$(g_{M/S},\omega_{M/S})$: with our normalization of $g_0$ and $g_{\infty}$
each fibre is the K\"ahler quotient of the corresponding fibre of $E_0\oplus
E_\infty$ by the diagonal $U(1)$ action at momentum level $z_0+z_\infty=2$;
then on this momentum level the function $z=z_0-1= 1-z_\infty$ descends to a
fibrewise momentum map $M\to [-1,1]$ for the quotient $U(1)$ action.

We extend $(g_{M/S},\omega_{M/S})$ to $TM$ by requiring that the horizontal
distribution of the induced connection on $M$ is degenerate. To obtain a
nondegenerate metric, we then set
\begin{equation*}
g_c = \sum_{a\in\cA} \frac{1+x_a z}{x_a} g_a + g_{M/S},
\qquad \omega_c=\sum_{a\in\cA} \frac{1+x_a z}{x_a} \omega_a + \omega_{M/S},
\end{equation*}
where the $(g_a, \omega_a)$ are pulled back from $S$; $g_c$ is then a positive
definite K\"ahler metric with respect to the canonical complex structure of
$M=P(\vE_0 \oplus \vE_{\infty})$ by the assumptions on the parameters $x_a$.
We refer to $(g_c,\omega_c)$ as the {\it canonical K\"ahler metric} on $M$ in
the given admissible K\"ahler class.

\begin{lemma} For any $0<|x_a|<1$ $(a\in\cA)$, the corresponding canonical
K\"ahler metric on $M$ is of the form~\eqref{metric} on $M^0$, where
$\Theta=\Theta_c$ and $\Theta_c(z)=1-z^2$.
\end{lemma}
\begin{proof} The inverse image in $E_0\oplus E_\infty$ of $M^0 = M \setminus
(z^{-1}(-1)\cup z^{-1}(1))$ may be viewed as an open subset of
$\cO(-1)_{E_0}\oplus\cO(-1)_{E_\infty}$. Then $(g_c,\omega_c)$ is the K\"ahler
quotient at momentum level $z_0+z_\infty=2$ of the metric
\begin{equation*}
\sum_{a} \frac{(1+x_a) z_0+(1-x_a)z_\infty}{2x_a} g_a + \frac{dz_0^2}{2z_0}
+ \frac{dz_\infty^2}{2z_\infty} + 2z_0\theta_0^2+2z_\infty \theta_\infty^2,
\end{equation*}
where $x_0=1$, $x_\infty=-1$, and $\theta_0$, $\theta_\infty$ are connection
$1$-forms for the $U(1)$-line bundles $\cO(-1)_{E_0}$, $\cO(-1)_{E_\infty}$
with $\theta_0(K_0)=1=\theta_\infty(K_\infty)$,
$d\theta_0=-\omega_0+\Omega_0$,
$d\theta_\infty=\omega_\infty+\Omega_\infty$.

If regard $M^0$ as an open subset of $\smash{\hat
M}=P(\cO\oplus\smash{\hat\cL})$, then the diagonal action is generated by
$K_0+K_\infty$, $\theta_\infty-\theta_0$ is basic and so induces a unitary
connection $\theta$ (with respect to the quotient $U(1)$-action) on
$\smash{\hat \cL}$ with $d\theta=\omega_{\hat S}$.  Substituting $z_0=1+z$ and
$z_\infty=1-z$ and performing the quotient yields~\eqref{metric} with
$\Theta=\Theta_c$. (On each fibre over $\smash{\hat S}$ this is the
realization of $\C P^1$ as a K\"ahler quotient of $\C^2$.)
\end{proof}

\begin{rem}\label{classes-adm} The existence of the canonical metric on $M$
shows there does exist a cohomology class $\Xi$ whose pullback to $\smash{\hat
M}$ is $\smash{\hat\Xi}$. $\Xi$ is then unique, and the admissible K\"ahler
classes form a family of dimension $\# \cA+1$.  If $b_2(S_a)=1$ for all $a$
and $b_1(S_a)\neq 0$ for at most one $a$, then every K\"ahler class on $M$ is
admissible.
\end{rem}

\subsection{Admissible metrics}\label{compactify}

Let $M=P(E_0\oplus E_\infty)\to S$ be an admissible bundle and $\Omega$ an
admissible K\"ahler class corresponding to parameters $x_a$. Then a K\"ahler
metric in $\Omega$ is said to be {\it admissible} if it has the
form~\eqref{metric} on $M^0$, up to scale, with respect to the given
projective unitary bundle structure on $M$ and local K\"ahler product
structure on $S$. According to Theorem~\ref{ACGTthm}, in order for a (scale of
a) metric of the form~\eqref{metric} on $M^0$ to define an admissible K\"ahler
metric on $M$, it is necessary that $\Theta$ is a smooth function on $[-1,1]$
satisfying~\eqref{positivity}--\eqref{boundary}. We now show that these
conditions are also sufficient and provide a parameterization of admissible
metrics.

We first note that any metric of the form~\eqref{metric}, where $\Theta$ is a
smooth function on $[-1,1]$ satisfying~\eqref{positivity}--\eqref{boundary},
defines a smooth metric $g$ on $M$ compatible with the same {\it symplectic
form} as the canonical K\"ahler metric $g_c$ in $\Omega$, provided that we
take $z$ to be the momentum map and $\theta$ the connection $1$-form of the
canonical K\"ahler metric; then, using~\eqref{boundary}, we find that $g-g_c$
is smooth on $M$, and $g$ is positive definite on $M$ since it is on $M^0$
by~\eqref{positivity} and $\omega$ is nondegenerate on $M$.
(See~\cite[\S1]{ACGT3} for details.)

With this point of view, the smooth functions $\Theta$ on $[-1,1]$
satisfying~\eqref{positivity}--\eqref{boundary} define a family of complex
structures on $M$. However, we claim that there is an $S^1$-equivariant
biholomorphism in the identity component of the diffeomorphism group between
any two such complex structures, so that $\Theta$ parameterizes K\"ahler
metrics compatible with the given (fixed) complex structure on $M$ whose
K\"ahler forms belong to a given admissible K\"ahler class $\Omega$.  This
claim holds essentially because it is true for toric complex structures on $\C
P^1$ (and for toric varieties in general), but for later use we need to make
explicit the transformation of $M$ that relates the complex and symplectic
points of view, following~\cite{guillemin,guan1,Do2}.

A key ingredient in this transformation is the notion of a {\it symplectic
potential} of an admissible K\"ahler metric defined by $\Theta(z)$, which is a
function $u(z)$ on $(-1,1)$ with $u''(z)=U(z)=1/\Theta(z)$. Then
\begin{equation*}
u_c(z) = \tfrac{1}{2}\bigl((1-z)\log (1-z) + (1+z) \log (1+z)-2\log 2\bigr)
\end{equation*}
is the unique symplectic potential for the canonical K\"ahler metric
$(g_c,J_c)$ given by $\Theta_c(z)= (1-z)(1+z)$, which satisfies $u_c(\pm
1)=0$.  We can extend this description to all admissible K\"ahler metrics
compatible with $\omega$, thanks to the following lemma, which is an easy
application of l'H\^opital's rule and Taylor's Theorem.
\begin{lemma}\label{parametrisation}
A smooth function $\Theta(z)=1/U(z)$ satisfies \eqref{boundary} if and only if
$U(z)-U_c(z)$ is smooth on $[-1,1]$. Then $U(z)/U_c(z)$ is positive and smooth
on $[-1,1]$.
\end{lemma}

On $M^0$ the symplectic potential $u(z)$ of an admissible K\"ahler metric is
closely related to a K\"ahler potential of $\omega$ with respect to $J$ by a
fibrewise Legendre transform (see \cite{ACG2,guillemin}) over $\smash{\hat
S}$. Indeed, if we put
\begin{equation*}
y = u'(z), \quad h(y)= -u(z) + yz,
\end{equation*}
then $d^c_J y=\theta$ and $dd^c_J h(y)=\omega-\sum_a \omega_a/x_a$ on
$(M^0,J)$\footnote{It follows that if $\pm H_a$ is a local K\"ahler potential
for $\pm\omega_a$ and $\tilde u = u(z) - \sum_a (1+x_a z) H_a/x_a$, then
$\tilde y = \partial \tilde u/\partial z$ is pluriharmonic and $\tilde
h=-\tilde u + \tilde yz$ is a local K\"ahler potential for $\omega$ on
$(M^0,J)$~\cite{ACG2}.}.  Let $y_c, h_c(y_c)$ denote the corresponding
quantities associated to $u_c$. There are local $1$-forms $\alpha$ on
$\smash{\hat S}$ such that $\theta=dt+\alpha$, where $t\colon M^0\to
\R/2\pi\Z$ is locally defined up to an additive constant on each fibre. Since
$\exp(y+it)$ and $\exp(y_c+it)$ give $\C\mult$-coordinates on the fibres,
there is a $U(1)$-equivariant fibre-preserving diffeomorphism $\Psi$ of $M^0$
over $\smash{\hat S}$ with
\begin{equation*}
\Psi^* y= y_c, \quad \Psi^*t = t, \quad \text{and hence}\quad \Psi^* J =J_c.
\end{equation*}
As $J_c$ and $J$ are integrable complex structures, $\Psi$ extends to a
$U(1)$-equivariant diffeomorphism of $M$ leaving fixed any point on $e_0 \cup
e_{\infty}$ (since it is fibre preserving).

Put $\tilde \omega := \Psi^*\omega$. Then $\tilde\omega$ is a K\"ahler form
on $(M,J_c)$ which (we claim) belongs to the same cohomology class $\Omega$ as
$\omega$. Indeed, on $M^0$ we have
\begin{equation*}
\tilde\omega - \omega = dd^c_{J_c} (h(y_c) - h_c(y_c))
\end{equation*}
since $dd^c_{J_c} h(y_c) = \Psi^* dd^c_{J} h(y)=\tilde\omega-\sum_a
\omega_a/x_a$, so the following implies the claim.
\begin{lemma} The function $h(y_c) - h_c(y_c)$ is smooth on $M$.
\end{lemma}
\begin{proof} Since $\Psi$ is a diffeomorphism with $\Psi^*y=y_c$, this holds
if and only if $h(y)-h_c(y)$ is smooth on $M$.  We already know that
$h(y)-h_c(y_c)= - (u(z)-u_c(z)) + z(u'(z)-u_c'(z))$ is smooth (by
Lemma~\ref{parametrisation}) so it suffices to show that $h_c(y)-h_c(y_c)$ is
smooth on $M$. However, knowing $u_c$ explicitly, we calculate
\begin{equation*}
h_c(y) - h_c(y_c) = - \frac{1}{2} \Big(\log\Big(\frac{1-{\tilde z}}{1-z}\Big)
+ \log \Big(\frac{1+{\tilde z}}{1+z} \Big)\Big), 
\end{equation*}
where $\tilde z := \Psi^* z$ is the momentum map of $\tilde\omega =
\Psi^*\omega$; since $\Psi$ is $S^1$-equivariant and fixes $e_0 \cup
e_{\infty}$, it follows that ${\tilde z}$, viewed as a function of $z$,
satisfies ${\tilde z}(\pm 1) = \pm 1$; moreover, since both ${\tilde z}$ and
$z$ are momentum maps of the same $U(1)$ action on $M$ (and are therefore
Morse-Bott functions with the same critical sets), we must have ${\tilde
z}'(\pm 1) \neq 0$.  Thus $h_c(y) - h_c(y_c)$ is smooth on $M$.
\end{proof}
Hence the moduli space ${\mathcal K}^{\rm adm}_{\omega}$ of admissible metrics
in $\Omega=[\omega]$ is identified with the space of smooth functions $\Theta$
on $[-1,1]$ satisfying \eqref{positivity}--\eqref{boundary} or equivalently
with $\{ u \in C^{0}([-1,1]) : u-u_c\in C^\infty([-1,1]), u(\pm 1)=0\text{ and
} u''>0 \text{ on } (-1,1) \}$.

\subsection{The isometry Lie algebra}\label{isometry}

For a compact K\"ahler manifold $(M,g)$, we denote by $\mathfrak{i}_0(M,g)$
the Lie algebra of all Killing vector fields with zeros. Since $M$ is compact
this is equivalently the Lie algebra of all hamiltonian Killing vector fields.

\begin{prop}\label{isometry algebra}
Let $g$ be an admissible metric on $M= P(E_0\oplus E_\infty)\stackrel{p}{\to}
S$ and equip $S$ and $\smash{\hat S}\to S$ with the metrics $g_S$, $g_{\hat
S}$ induced by $\sum_a g_a/x_a$ on $\prod_a S_a$. Let $\mathfrak z(K,g)$ be
the centralizer in $\mathfrak{i}_0(M,g)$ of the Killing vector field
$K=J\grad_g z$.

Then $\mathfrak z(K,g)$ is the direct sum of $\mathfrak{i}_0(\smash{\hat
S},g_{\hat S})$ and the span of $K$ in such a way that $p_*\colon
\mathfrak{i}_0(M,g)\to \mathfrak{i}_0(S,g_S)$ is induced by the natural
surjection $\mathfrak{i}_0(\smash{\hat S},g_{\hat S})\to
\mathfrak{i}_0(S,g_S)$.
\end{prop}
\begin{proof} Let $X$ be a holomorphic vector field on $\smash{\hat S}$
which is hamiltonian with respect to $\omega_h:=\sum_a\omega_a/x_a$; then the
projection $X_a$ of $X$ onto the distribution ${\mathcal H}_a$ (induced by
$TS_a$ on the universal cover $\prod_a S_a$ of $\smash{\hat S}$) is a Killing
vector field with zeros, so $\iota^{\vphantom{x}}_{X_a}\omega_h=-df_a$ for
some function $f_a$ (with integral zero). Thus $\sum_a f_ax_a$ is a
hamiltonian for $X$ with respect to the symplectic form $\omega_{\hat
S}=\sum_a \omega_a$: since this is the curvature $d\theta$ of the connection
on $M^0$, $X$ lifts to a holomorphic vector field $\tilde X =X_H+(\sum_a
f_ax_a) K$ on $M^0$, which is hamiltonian with potential $\sum_a(1+ x_a z)f_a$
and commutes with $K$. Here $X_H$ is the horizontal lift to $M^0$ with respect
to $\theta$.  $\tilde X$ and its potential extend to $M$ since $M\setminus
M^0$ has codimension $\geq 2$ and $\tilde X$ has zeros.

Conversely any element of $\mathfrak z(K,g)$ pulls back to a holomorphic
vector field $V$ on $\hat M$.  The projection of $V$ to the normal bundle
$\hat p^*T\smash{\hat S}$ of $\hat p\colon\smash{\hat M}\to \smash{\hat S}$ is
holomorphic hence constant on the $\C P^1$ fibres by Liouville's Theorem (the
normal bundle is trivial on each fibre), so $V$ is projectable; it maps to
zero iff it comes from a multiple of $K$. This gives a projection to
$\mathfrak{i}_0(\hat S,g_{\hat S})$ splitting the inclusion just defined.
\end{proof}

\section{Admissible extremal K\"ahler metrics}\label{main}

\subsection{Automorphisms and the Matsushima--Lichnerowicz obstruction}
\label{matsushima-lichnerowicz}

On any compact K\"ahler manifold $(M,g)$ the Lie algebra $\mathfrak{h}(M)$ of
holomorphic vector fields lies in an exact sequence:
\begin{equation*}
0\to \mathfrak{h}_0(M) \to \mathfrak{h}(M) \to H^1(M,\R)^*
\end{equation*}
where $\mathfrak{h}_0(M)$ is the ideal of holomorphic vector fields with
zeros, which is the Lie algebra of the reduced automorphism group
$H_0(M)\subset\Aut_0(M)$, the connected component of the kernel of the
Albanese map ${\Aut}_0(M)\to H^1(M,\R)^*/H_1(M,\Z)$.  The
Matsushima--Lichnerowicz Theorem~\cite{lichne,matsushima} says that if $g$ is
CSC, $\mathfrak{h}_0(M)$ is the complexification of the Lie algebra
$\mathfrak{i}_0(M,g)$ of hamiltonian Killing vector fields and
$\mathfrak{h}(M)=\mathfrak a(M)\oplus\mathfrak h_0(M)$ where $\mathfrak{a}(M)$
is the central subalgebra of parallel vector fields: thus $\mathfrak{h}(M)$ is
reductive.  This condition on $\mathfrak{h}(M)$ is often called the {\it
Matsushima--Lichnerowicz obstruction} to the existence of CSC K\"ahler metrics
on $(M,J)$.

Let us consider the special case of geometrically ruled complex manifolds
$M=P(E)$ where $p\colon P(E)\to S$ and $E$ is a holomorphic vector bundle of
rank $r+1$ over a compact K\"ahler $2d$-manifold $(S,h,\omega_h)$. Thus, $M$
is a k\"ahlerian $2m$-manifold where $m=r+d$: for instance, $[\cO(1)_E] + k
[p^*\omega_h]$ is a K\"ahler class for $k\gg 1$, where $\cO(-1)_E$ is the
(fibrewise) tautological line bundle of $P(E)$.

The projection of any holomorphic vector field $V\in\mathfrak{h}(M)$ to the
normal bundle $p^*TS$ is constant on each fibre, so $V$ descends to a
holomorphic vector field $p_*V\in\mathfrak h(S)$. Since
$p_*[V_1,V_2]=[p_*V_1,p_*V_2]$, we have an exact sequence of Lie algebras
\begin{equation*}
0 \to  \mathfrak{h}_S(M) \to \mathfrak{h}_0(M)\to \mathfrak{h}_0(S),
\end{equation*}
where $\mathfrak{h}_S(M)$ is the subspace of $\mathfrak{h}(M)$ of holomorphic
vector fields tangent to the fibres of $p$ (which have zeros).  Obviously
$\mathfrak{h}_{S}(M)=H^0(S, \mathfrak{sl}(E))$ is the Lie algebra of
holomorphic vector fields preserving the ${\C}P^{r}$-fibres of $p\colon
P(E)\to S$. Since an ideal in a reductive Lie algebra is reductive, we obtain
the following weaker (but often more useful) version of the
Matsushima--Lichnerowicz obstruction.

\begin{prop}\label{relative} Let $M= P(E) \to S$ be a geometrically ruled
complex manifold which admits a CSC K\"ahler metric. Then $\mathfrak{h}_S(M)$
must be reductive.
\end{prop}

The following elementary result yields a simple application of this criterion.

\begin{prop}\label{reductive} Let $M=P(\vE_0\oplus \vE_1\oplus \cdots \oplus
\vE_\ell)$, where $\vE_j$ are holomorphic vector bundles on a compact
k\"ahlerian manifold $S$ such that $H^0(S,\mathfrak{gl}(\vE_j))$ is reductive
and $H^0(S,\Hom(\vE_{j},\vE_i))=0$ for all $i<j$.  Then ${\mathfrak h}_S(M)$
is reductive iff $H^0(S,\Hom(\vE_{i},\vE_j))=0$ for all $i<j$.
\end{prop}
\begin{proof} Any element of the Lie algebra $\mathfrak{h}_{S}(M)
=H^0(S,\mathfrak{sl}(E))$ may be represented as an $(\ell+1)\times (\ell+1)$
matrix $(a_{ij})$ with $a_{ij}\in H^0(S,\Hom(\vE_{i},\vE_j))$. By assumption,
this matrix is upper-triangular.  The strictly upper-triangular matrices form
a nilpotent ideal $\mathfrak n$ and if this is zero, $\mathfrak{h}_{S}(M)$ is
clearly reductive.  Conversely, taking commutators with elements of the form
$\Id_{E_i}/\rk E_i-\Id_{E_j}/\rk E_j\in H^0(S,\mathfrak{sl}(E))$,
we see that $\mathfrak n\subseteq [\mathfrak{h}_S(M), \mathfrak{h}_S(M)]$.
Hence if $\mathfrak n\neq 0$, $[\mathfrak{h}_S(M), \mathfrak{h}_S(M)]$ is not
semisimple, i.e., $\mathfrak{h}_S(M)$ is not reductive.
\end{proof}

\begin{cor}\label{rscase} Let $M=P(\vE)$ where $\vE=\cL_0 \oplus \cL_1 \oplus
\cdots \oplus \cL_r$ is a direct sum of holomorphic line bundles over a
Riemann surface $\Sigma$ of genus $\mathbf g$. If $\mathbf g\geq 2$ and
$\deg\cL_i-\deg\cL_j> {\mathbf g} -1$ for some $0\leq i,j\leq r$ then $M$
admits no CSC K\"ahler metric. If ${\mathbf g}\leq 1$, then $M$ admits a CSC
K\"ahler metric if and only if $\deg\cL_i = \deg\cL_j$ for all $i,j$
\textup(i.e., $\cL_0 \oplus \cL_1 \oplus \cdots \oplus \cL_r$ is a polystable
vector bundle\textup).
\end{cor}
\begin{proof} We can assume without loss that $E=\lE_0\otimes
\C^{r_0} \oplus \lE_1 \otimes \C^{r_1} \oplus \cdots \oplus \lE_\ell\otimes
\C^{r_{\ell}}$ with $\deg\lE_i \leq \deg\lE_j$ and $\lE_i \not\cong \lE_j$ for
$i<j$.  The Kodaira vanishing theorem then implies $H^0(\Sigma, \lE_{\smash
j}^{-1}\otimes \lE_i)=0$ for any $i<j$, and we may apply
Proposition~\ref{reductive}.  By assumption, $\deg\lE_\ell- \deg\lE_0 >
{\max}(0, {\mathbf g} -1)$, and so $\dim H^0(\Sigma,
\lE_0^{-1}\otimes\lE_\ell) >0$ by Riemann--Roch. Hence ${\mathfrak
h}_\Sigma(M)$ is not reductive, and there is no CSC K\"ahler metric on
$M$. The converse when $\mathbf g\leq 1$ follows from
Narasimhan--Seshadri~\cite{ns}.
\end{proof}

\begin{rem}\label{flat bundles}
The assumptions of Proposition~\ref{reductive} hold if $\vE_1,\ldots \vE_\ell$
are projectively-flat hermitian vector bundles such that the slopes
$\mu(\vE_j):=\cb_1(\vE_j)\cup\Omega^{d-1}$, with respect to some K\"ahler
class $\Omega$ on $S$ ($\dim S=2d$), satisfy $\mu(\vE_i)<\mu(\vE_j)$ for
$i<j$.  Indeed in this case $\mathfrak{gl}(\vE_j)$ is a flat hermitian bundle
and $H^0(S,\mathfrak{gl}(\vE_j))$ is the space of parallel sections of
$\mathfrak{gl}(\vE_j)$\footnote{This is a standard Bochner argument, as in
\cite{kob}. Alternatively, note that the pullback of $\mathfrak{gl}(\vE_j)$ to
the universal cover of $S$ is trivialized by parallel sections, and apply the
open mapping theorem: the pullback of a holomorphic section of
$\mathfrak{gl}(\vE_j)$ has closed bounded image in this trivialization.},
which is a complexification of the space of parallel sections of $\mathfrak
u(\vE_j)$, hence a reductive Lie algebra.  The slope condition then ensures
$H^0(S,\Hom(\vE_{j},\vE_i))=0$ for all $i<j$ by a theorem of
Kobayashi~\cite{kob}.
\end{rem}

In general, the condition that $\mathfrak{h}_S(M)$ is reductive does not imply
$\mathfrak h_0(M)$ is. However it does if $p_*\colon\mathfrak{h}_0(M)\to
\mathfrak{h}_0(S)$ is surjective and $\mathfrak{h}_0(S)$ is reductive.  This
obviously holds if $\mathfrak{h}_0(S)=0$. It also holds if $(S,g_S)$ is CSC
and there is a metric $g$ on $M$ such that $p_*$ is a surjection from
$\mathfrak{i}_0(M,g)$ to $\mathfrak{i}_0(S,g_S)$. This is true for admissible
bundles by Proposition~\ref{isometry algebra}.

\begin{prop}\label{maximal} Let $\Omega$ be an admissible K\"ahler class on
$M= P(\vE_0\oplus\vE_\infty)\to S$ where the local product metric $g_S$ on $S$
is CSC.  Then the admissible metrics in $\Omega$ are invariant under a common
maximal compact connected subgroup of $H_0(M)$, and $\mathfrak h_0(M)$ is
reductive iff $H^0(S,\Hom(E_0, E_\infty))=0=H^0(S,\Hom(E_\infty,E_0))$. This
latter condition holds if $\cb_1(\vE_\infty)-\cb_1(\vE_0)$ is strictly
indefinite.
\end{prop}
\begin{proof} Let $g$ be an admissible K\"ahler metric on $M$. We know that
$\mathfrak h(S)$ is reductive and $\mathfrak i_0(S,g_S)$ is a maximal compact
subalgebra of $\mathfrak h_0(S)$, and by Proposition~\ref{isometry algebra},
both $p_*\colon\mathfrak{h}_0(M)\to \mathfrak{h}_0(S)$ and
$p_*\colon\mathfrak{i}_0(M,g)\to\mathfrak{i}_0(S,g_S)$ are surjective.

To show that the Lie algebra $\mathfrak i_0(M,g)$ is a maximal compact
subalgebra of $\mathfrak h_0(M)$, it therefore suffices to show that
$\mathfrak{i}_0(M,g)\cap\mathfrak h_S(M)$ is a maximal compact subalgebra of
$\mathfrak h_S(M)$. Since $\cb_1(\vE_\infty)-\cb_1(\vE_0)=\omega_S$, we can
certainly arrange that $\mu(\vE_\infty)-\mu(\vE_0)$ is nonzero by the choice
of a K\"ahler class on $S$. Then, by Remark~\ref{flat bundles}, we are under
the hypotheses of Proposition~\ref{reductive}, and, as in its proof, we have
that $\mathfrak h_S(M)$ is the direct sum of the reductive centralizer of $K$
and a nilpotent ideal $\mathfrak n$ in $[\mathfrak h_S(M), \mathfrak h_S(M)]$.
The result now follows easily from Proposition~\ref{isometry algebra} and
Remark~\ref{flat bundles}.

As noted above, $\mathfrak h_0(M)$ is reductive iff $\mathfrak{h}_S(M)$ is.
By Proposition~\ref{reductive}, the latter happens iff both $H^0(S,\Hom(E_0,
E_\infty))=0=H^0(S,\Hom(E_\infty,E_0))$. This indeed holds if
$\cb_1(\vE_\infty)-\cb_1(\vE_0)$ is strictly indefinite (by the vanishing
theorem of Kobayashi~\cite{kob} as in~Remark 3), since we can then choose
K\"ahler classes on $S$ such that the corresponding slopes have
$\mu(\vE_\infty)-\mu(\vE_0)$ with either sign.
\end{proof}
Since $S$ is a local K\"ahler product, it is CSC iff the factors
$S_a$ $(a\in\cA)$ in the universal cover are CSC.

\subsection{The Futaki invariant and extremal vector field}
\label{futaki-invariant}

On a compact K\"ahler $2m$-manifold $(M,J,g,\omega)$, recall that the
(normalized) {\it Futaki invariant} of a real holomorphic vector field with
zeros $V = J \grad_g f + \grad_g h$ is defined by
\begin{equation*}
\mathfrak{F}_{\omega} (V) = \biggl({\int_M \mu_g}\,
{\int_M (f+ih){\Scal_g} \mu_g}
- {\int_M {\Scal_g} \mu_g}\,{\int_M (f+ih)\mu_g}
\biggr)\big/\Vol(M)^2,
\end{equation*}
where $\mu_g=\omega^m/m!$ is the volume form of $g$. Futaki~\cite{futaki}
showed that this complex number is independent of the choice of metric in the
K\"ahler class $\Omega=[\omega]$, and that the map
$\mathfrak{F}_{\Omega}\colon \mathfrak{h}_0(M) \to \C$ is a character on
$\mathfrak h_0(M)$.  $\mathfrak{F}_{\Omega}$ is closely related to the
{\it Futaki--Mabuchi extremal vector field}
$K_{\Omega}:=J\grad_g\pr_g\Scal_g$ of $(M,J,\Omega,G)$ where $G$ is a maximal
compact connected subgroup of $H_0(M)$ and $\pr_g$ is the $L_2$-projection
onto the space of Killing potentials with respect to any $G$-invariant metric
$g$ in $\Omega$: Futaki and Mabuchi~\cite{futaki-mabuchi} showed that
$K_{\Omega}$ is independent of this choice. Clearly $\mathfrak{F}_{\Omega}$
and $K_{\Omega}$ vanish if $\Omega$ contains a CSC metric.
Calabi~\cite{calabi2} showed that if $\mathfrak{F}_{\Omega}$ vanishes then any
extremal K\"ahler metric in $\Omega$ is a CSC metric, but the vanishing of
$\mathfrak{F}_{\Omega}$ does not suffice in general for the existence of a CSC
metric in $\Omega$.

Let $\Omega$ be an admissible K\"ahler class on $M=P(\vE_0 \oplus\vE_\infty)
\to S$ and suppose in addition that for $a\in\cA$, $\pm g_a$ is a CSC K\"ahler
metric with scalar curvature $\Scal_{\pm g_a}=\pm 2d_as_a$. Let
$\Mpc(t)=\prod_a (1+x_at)^{d_a}$ and define $\alpha_r = \int_{-1}^{1}
\Mpc(t)t^r dt$ and
\begin{equation}\label{betas}
\beta_r:= \Mpc(1)+(-1)^r\Mpc(-1)
+\smash[t]{\int_{-1}^1 \biggl(\sum_{a} \frac{d_a s_a x_a}{1+x_a t}\biggr)}
\vphantom{\Big|}\Mpc(t) t^r dt.
\end{equation}
We now compute the Futaki invariant ${\mathfrak F}_{\Omega}(K)$ of $K=J\grad_g
z$ and show that $K_{\Omega}$ is essentially ${\mathfrak F}_{\Omega}(K) K$,
where $G$ the maximal compact connected subgroup of $H_0(M)$ of
Proposition~\ref{maximal} preserving admissible K\"ahler metrics in $\Omega$.
${\mathfrak F}_{\Omega}(K)$ will reappear in the next paragraph as the leading
coefficient of a polynomial associated with $\Omega$.

\begin{prop}\label{futaki-ext} Suppose $M$ is admissible over a CSC base and
$\Omega$ is an admissible K\"ahler class with admissible metric $g$. Then
$\mathfrak{F}_{\Omega}(K) = 2(\alpha_0\beta_1-\alpha_1\beta_0) /\alpha_0^2$.
Also the $L_2$-projection of $\Scal_g$ orthogonal to the space of Killing
potentials is
\begin{equation}\label{prScal}
\Scal_g+A z+B
\end{equation}
where $A$ and $B$ are given by
\begin{equation}\label{system}\begin{split}
A \alpha_1 + B \alpha_0 &= -2\beta_0\\
A \alpha_2 + B \alpha_1 &= -2\beta_1.
\end{split}\end{equation}
\textup(Since $\alpha_0\alpha_2>\alpha_1^2$, this system has a unique solution
for $A,B$.\textup) In particular the extremal vector field of $(\Omega,G)$ is
$K_{\Omega}=-A K =2(\alpha_0\beta_1-\alpha_1\beta_0)K
/(\alpha_0\alpha_2-\alpha_1^2)$.
\end{prop}
\begin{proof} We may rescale $\Omega$ so that an admissible metric
$(g,\omega)$ in $\Omega$ is of the form~\eqref{metric}. We then have
\begin{equation*}
\mu_g=\frac{\omega^m}{m!} = \Mpc(z)
\biggl(\bigwedge_{\smash a} \frac{(\omega_a/x_a)^{d_a}}{d_a!}\biggr)
\wedge dz\wedge \theta,
\end{equation*}
where $m = 1+\sum_a d_a$ is the complex dimension of $M$. Thus
\begin{align*}
\int_M \mu_g &
= 2\pi \Vol\bigl(S,{\textstyle \prod_a\frac{\omega_a}{x_a}}\bigr) \alpha_0
\,\bigl(=\Vol(M)\bigr),\\
\int_M z \mu_g &
= 2\pi \Vol\bigl(S,{\textstyle \prod_a\frac{\omega_a}{x_a}}\bigr) \alpha_1
=\Vol(M)\alpha_1/\alpha_0,
\end{align*}
where $\Vol\bigl(S,{\textstyle \prod_a\frac{\omega_a}{x_a}}\bigr)= \prod_a
\Vol\bigl(S_a,{\textstyle \frac{\omega_a}{x_a}}\bigr)$ in the case when $S$ is
a global product.

The scalar curvature of $(M,g)$ is given by
\begin{equation}\label{scal-form}
{\Scal_g}= \sum_a \frac{2 d_a s_a x_a }{1+x_az} - \frac{F''(z)}{\Mpc(z)}
\end{equation}
where $F(z)=\Theta(z)\Mpc(z)$ (see e.g.~\cite[(79)]{ACG2}). We thus calculate
\begin{align*}
\int_M &z {\Scal_g} \mu_g = 2\pi \Vol\bigl(S,{\textstyle
\prod_a\frac{\omega_a}{x_a}}\bigr) \int_{-1}^{1}\biggl(
\Big(\sum_a \frac{2 d_a s_a x_a}{1+x_a z}\Bigr)\Mpc(z) - F''(z)\biggr)z dz\\
&=2\pi \Vol\bigl(S,{\textstyle \prod_a\frac{\omega_a}{x_a}}\bigr)
\biggl( \int_{-1}^{1}
\Bigl( \sum_a \frac{2d_a s_a x_a }{1+x_a z}\Bigr)\Mpc(z)z\, dz -
\Bigl[zF'(z) - F(z)\Bigr]^1_{-1}\biggr)\\
&= 2 \Vol(M)\beta_1/\alpha_0,
\end{align*}
where we integrate by parts, then impose the boundary conditions
\eqref{boundary}.  Similarly,
\begin{equation*}
\int_M \Scal_g \mu_g = 2\Vol(M)\beta_0/\alpha_0
\end{equation*}
and the first claim follows.

For the second claim note that the above integral formulae imply $\Scal_g+A
z+B$ is orthogonal to the Killing potentials $1,z$ if and only
if~\eqref{system} holds. By the form of $\Scal_g$, the fact that the $s_a$ are
constant, and Proposition~\ref{isometry algebra}, the result follows.
\end{proof}
Note that the above expression for $\mathfrak{F}_{\Omega}(K)$ is manifestly
independent of the choice of a smooth function $\Theta(z)$ satisfying
\eqref{boundary}, as it should be according to the general
theory~\cite{futaki}. Indeed, as we have already discussed in
\S\ref{compactify}, these smooth functions $\Theta(z)$ define K\"ahler metrics
within the same K\"ahler class.

\subsection{K-energy and the extremal polynomial}
\label{K-energy-ext}

Given a complex $2m$-manifold $(M,J)$, a maximal compact connected subgroup
$G$ of $H_0(M)$ and a K\"ahler class $\Omega$, denote by ${\mathcal
M}_{\Omega}$ the infinite dimensional Fr\'echet space of K\"ahler metrics in
$\Omega$ and let ${\mathcal M}^G_{\Omega}$ be the subspace of $G$-invariant
K\"ahler metrics in $\Omega$. Following Guan~\cite{guan1} and
Simanca~\cite{simanca}, consider the map
\begin{equation*}
g \mapsto \pr_g^\perp\Scal_g \mu_g,
\end{equation*}
where $\pr_g^\perp$ is the $L_2$-projection orthogonal to the space of Killing
potentials. This can be viewed (by integration) as a 1-form $\sigma$ on
${\mathcal M}_{\Omega}^G$, which turns out to be closed. Therefore for any
$\omega_0\in \Omega$, there exists a unique functional ${E}
_{\omega_0}^G\colon {\mathcal M}^G_{\Omega} \to \R$ with
\begin{equation*}
d {E} _{\omega_0}^G= - {\sigma},
\end{equation*}
${E} _{\omega_0}^G(\omega_0)=0$. Note that changing the base point $\omega_0
\in {\mathcal M}_{\Omega}$ would change $E _{\omega_0}^G$ by an additive
constant. We refer to $E_{\omega_0}^G$ as the (modified) {\it K-energy}: it
agrees with the Mabuchi K-energy~\cite{mabuchi-0} when $G$ is trivial.

By definition, it is clear that the critical points of ${E}^G _{\omega_0}$ are
exactly the extremal K\"ahler metrics in ${\mathcal M}_{\Omega}^G$, since
$\sigma=0$ means that $\Scal_g$ is a Killing potential. Note that by the
Calabi Theorem~\cite{calabi1}, any extremal K\"ahler metric $g\in {\mathcal
M}_{\Omega}$ belongs to ${\mathcal M}_{\Omega}^G$ with $G= {\rm Isom}_0(M,g)
\cap H_0(M)$.

Building on earlier work by Bando--Mabuchi~\cite{Bando-Mabuchi},
Chen~\cite{chen1}, Donaldson~\cite{Do3} and others, Chen and Tian have
established the following uniqueness result and necessary condition for
existence of an extremal K\"ahler metric.

\begin{thm}\cite{CT,CT2} \label{K-energy-thm} Extremal K\"ahler metrics in
${\mathcal M}_\Omega$ are unique up to automorphism and any extremal K\"ahler
metric in ${\mathcal M}_\Omega^G$ realizes the absolute minimum of
${E}_{\omega_0}^G$ \textup(for any $\omega_0 \in{\mathcal M}_{\Omega}^G
$\textup). In particular, if ${\mathcal M}_{\Omega}^G$ contains an extremal
K\"ahler metric, then ${E}_{\omega_0}^G$ is bounded from below.
\end{thm}

Now let $M$ be an admissible projective bundle over a CSC base as in the
previous paragraph.  We want to obtain a formula for the K-energy as a
functional acting on ${\mathcal K}^{\rm adm}_{\omega}$, where $\omega$ is
fixed, so we need to use the description given in \S\ref{compactify} which
shows how ${\mathcal K}^{\rm adm}_{\omega}$ is embedded into ${\mathcal
M}_{\Omega}^G$, in which the complex structure is fixed.

This description shows that if $u_t(z)$ is a path of symplectic potentials in
${\mathcal K}^{\rm adm}_{\omega}$, then the smooth functions $h_t(y_c) -
h_c(y_c)$ define a path $\omega + dd^c_{J_c} (h_t(y_c) - h_c(y_c))$ in
${\mathcal M}_{\Omega}$, where $h_t$ are introduced by
\begin{equation*}
h_t(y_t(z)) = -u_t(z) + y_t z, \quad y_t = u_t'(z)
\end{equation*}
so that with $z=y_t^{-1}(y_c)$ we have
\begin{equation*}
h_t(y_c) = -u_t(y_t^{-1}(y_c)) + y_t^{-1}(y_c) y_c.
\end{equation*}
Differentiating with respect to $t$, we get for the corresponding vector
fields ${\dot u} \in T_{g}({\mathcal K}^{\rm adm}_{\omega})$ and ${\dot h}
\in T_{\omega} ({\mathcal M}_{\Omega}^G)$, the relation
(cf.~\cite{guan1,Do2}):
\begin{equation*}
\dot h = - \dot u.
\end{equation*}
Hence we obtain the following symplectic version of the (modified) K-energy.
\begin{lemma} \label{symplectic-setting} The K-energy ${E}_{\omega}^G$,
restricted to the space of admissible K\"ahler metrics in $\Omega$ and viewed
as a function on the space of symplectic potentials, is determined uniquely up
to an additive constant by the formula
\begin{equation*}
(d {E}_{\omega}^G)_{g}[\dot u] = \int_M (\pr_g^\perp \Scal_g) \dot u
\, \mu_g,
\end{equation*}
where $\pr_g^\perp$ denotes the $L_2$-projection orthogonal to the space of
Killing potentials.
\end{lemma}

Consider an admissible metric $g$ in $\Omega$ corresponding to the function
$\Theta(z)=F(z)/\Mpc(z)$.  Since the base $S$ is CSC we have, by
Proposition~\ref{futaki-ext}, $\pr_g^\perp\Scal_g =\Scal_g + A z + B$, with
$A$ and $B$ given by \eqref{system} and $\Scal_g$ by~\eqref{scal-form}.
\begin{lemma} There is a unique smooth function $F_\Omega$ on $[-1,1]$ with
\begin{equation}\label{ODE}
F''_\Omega(z)= \biggl(A z+ B+ \sum_{a} \frac{2 d_a s_a x_a}{1+x_a z}\biggr)
\Mpc(z)
\end{equation}
and $F_\Omega(\pm1)=0$. $F_\Omega$ satisfies~\eqref{BCs} and is a polynomial
of degree $\leq m+2$, the coefficient of $z^{m+2}$ being a nonzero multiple of
$A$.
\end{lemma}
\begin{proof} There is clearly a unique solution to~\eqref{ODE} with
$F_\Omega(\pm1)=0$.  One easily checks, using~\eqref{system} that the
solution is
\begin{equation*}
(1+z)\Mpc(1)+(1-z)\Mpc(-1)+ \int_{-1}^1 \biggl(\frac12(A t+ B)+ \sum_{a}
\frac{ d_a s_a x_a}{1+x_a t}  \biggr) \Mpc(t) |z-t| dt.
\end{equation*}
The derivative of this function is
\begin{equation*}
\Mpc(1)-\Mpc(-1)+\int_{-1}^1 \biggl(\frac12(A t+ B)+ \sum_{a}
\frac{ d_a s_a x_a}{1+x_a t}  \biggr) \Mpc(t) \mathrm{sign}(z-t) dt,
\end{equation*}
which gives the formulae for $F_\Omega'(\pm 1)$ in~\eqref{BCs}, using the
first equation of~\eqref{system}.
\end{proof}
The motivation for this lemma is that now $\Scal_g + A z + B =
(F_\Omega''(z)-F''(z))/\Mpc(z)$. Furthermore, $F$ and $F_\Omega$ satisfy the
same boundary conditions~\eqref{BCs}.
\begin{prop}\label{K-energy} Let $\Omega$ be an admissible K\"ahler
class on an admissible bundle over a CSC base. Then the K-energy restricted to
the space of admissible K\"ahler metrics ${\mathcal K}^{\rm adm}_\omega$ is
\textup(up to an additive constant\textup) a positive multiple of the
functional
\begin{equation*}
{\mathcal E}_{g_c}: u(z) \mapsto
\int_{-1}^1 F_\Omega(z) (u''(z)-u_c''(z)) dz 
-\int_{-1}^1 p_c(z) \log \Big(\frac{u''(z)}{u_c''(z)}\Big) dz.
\end{equation*}
\end{prop}
\begin{proof} ${\mathcal E}_{g_c}$ is well-defined by
Lemma~\ref{parametrisation} and its gradient is
\begin{align*}
(d {\mathcal E}_{g_c})_{g}[\dot u] &=\int_{-1}^1 F_\Omega(z) \dot
u''(z) dz -\int_{-1}^1 p_c(z) \frac{\dot u''(z)}{u''(z)} dz
=\int_{-1}^1(F_\Omega(z) -F(z)) \dot u''(z) dz.
\end{align*}
Integrating twice by parts, using the fact that $F$ and $F_\Omega$ both
satisfy~\eqref{BCs}, and multiplying by $2\pi {\rm Vol}(S,
\prod_a{\omega_a}/{x_a})$, we obtain $\int_M (Scal_g + A z + B)\dot u
\mu_g$.
\end{proof}
\begin{rem}\label{absolute Mabuchi}
It is worth noticing that (by~\eqref{boundary}) the integral
\begin{equation*}
\int_{-1}^{1} \bigl(F_{\Omega}(z) u''(z) -  p_c(z)\log u''(z)\bigr) dz
\end{equation*}
exists for any admissible K\"ahler metric in $\Omega$ with symplectic
potential $u(z)$, giving a definition of K-energy on ${\mathcal
K}_{\omega}^{\rm adm}$ which is independent of a choice of reference metric.
\end{rem} 
\begin{cor}\label{F-nonneg} If there is an extremal K\"ahler metric in
$\Omega$, then $F_\Omega\geq 0$ on $[-1,1]$.
\end{cor}
\begin{proof}
If there is an extremal K\"ahler metric in $\Omega$, then by
Theorem~\ref{K-energy-thm}~\cite{CT,CT2}, the K-energy is bounded from
below. We now apply an argument from~\cite{Do2}: take any nonnegative smooth
function $f(z)$ with ${\rm supp}(f)\subset (-1,1)$ and consider the sequence
$u_k(z)$ with $u_k''(z) = u_c''(z) + kf(z)$ of symplectic potentials
(cf.~Lemma~\ref{parametrisation}) for admissible K\"ahler metrics. We
therefore get
\begin{equation*}
{\mathcal E}_{g_c} (u_k) = -\int_{-1}^{1} p_c(z)\log
\Bigl(1+k\frac{f(z)}{u_c''(z)}\Bigr)dz + k\int_{-1}^{1} F_\Omega(z)f(z) dz.
\end{equation*}
This will tend to $-\infty$ if $\int_{-1}^{1} F_\Omega(z)f(z) dz <0$ for
some $f$.
\end{proof}

In the next paragraph, where we complete the proof of Theorem~\ref{main-thm},
we shall show that positivity of $F_{\Omega}$ on $(-1,1)$ is a necessary and
sufficient condition for the existence of an extremal K\"ahler metric in
$\Omega$.

\begin{defn} Let $\Omega$ be an admissible K\"ahler class on $M$. Then
the polynomial $F_\Omega$ constructed above will be called the {\it extremal
polynomial} of $\Omega$.
\end{defn}

\subsection{A characterization of extremal admissible K\"ahler classes}
\label{char}

In this paragraph we prove Theorem~\ref{main-thm} in three steps. First, if
the extremal polynomial $F_\Omega$ of an admissible K\"ahler class $\Omega$ is
positive on $(-1,1)$, we construct an admissible extremal K\"ahler metric in
$\Omega$ by adapting an argument of Guan and Hwang
(cf.~\cite{guan0,hwang,hwang-singer}): we discuss their work further in the
next section. Second, we extend the continuity argument of~\cite{christina1}
to prove the existence of admissible extremal K\"ahler metrics for $\Omega$
sufficiently small. Third, we use Corollary~\ref{F-nonneg}, the uniqueness
result of Chen--Tian~\cite{CT,CT2} and an argument from~\cite{christina1} to
show that an extremal K\"ahler metric in $\Omega$ is admissible up to
automorphism. Hence we deduce that the existence of an extremal K\"ahler
metric in $\Omega$ implies that $F_\Omega$ is positive on $(-1,1)$.

We begin with the construction. By Proposition~\ref{equations}, an admissible
metric~\eqref{metric} is extremal exactly when for each $a\in\smash{\hat\cA}$,
$\pm g_a$ is a CSC K\"ahler metric with $\Scal_{\pm g_a}=\pm 2d_as_a$ and
\eqref{extremal-deriv}--\eqref{extremal-values} hold for a polynomial $P$ of
degree $\leq N+1$, where $N=\#\smash{\hat\cA}$.  The metric $g$ is CSC iff $P$
has degree $\leq N$.

We have seen that the boundary conditions~\eqref{boundary} imply~\eqref{BCs}
and the converse clearly holds if $\Mpc(\pm1)\neq0$ (i.e., $d_0=0=d_\infty$).
However, if $g$ is extremal,
then~\eqref{extremal-deriv}--\eqref{extremal-values} imply that
$F''(z)=\Mpc'(z)\Upsilon(z)$ with $\Upsilon(-1)=2(d_0+1)$ if $d_0>0$ and
$\Upsilon(1)= -2(d_{\infty}+1)$ if $d_\infty>0$ (because of the normalization
of the Fubini--Study metrics on $S_0$ and $S_\infty$). Hence, by~\eqref{BCs},
$F'(z)=\Mpc(z)\Psi(z)$ with $\Psi(-1)= 2(d_{0}+1)$ and $\Psi(1)=
-2(d_{\infty}+1)$, and $\Theta(\pm 1)=0$. Now by l'H\^opital's rule,
$\Theta'(\pm 1)=\mp 2$. Hence for extremal K\"ahler metrics, the boundary
conditions~\eqref{boundary} are equivalent to~\eqref{BCs}.

In summary, to obtain a globally defined admissible extremal metric on a
projective bundle $P(E_0 \oplus E_{\infty}) \to S$, we need, for CSC K\"ahler
metrics $(\pm g_a,\pm \omega_a)$ satisfying $\cb_1(E_{\infty}) - \cb_1(E_0) =
\sum_a[\omega_a/2\pi]$, to solve \eqref{extremal-deriv} and
\eqref{extremal-values} for a polynomial $F$ (of degree $\leq m+2$) which
satisfies~\eqref{BCs} and is positive on $(-1,1)$.

For an admissible K\"ahler class $\Omega$ on $M$, we claim that
\eqref{extremal-deriv}--\eqref{extremal-values} and the boundary
conditions~\eqref{BCs} have a unique solution for $F$, given by the extremal
polynomial $F_\Omega$.

\begin{prop}\label{extremal-poly}
Let $M = P(\vE_0 \oplus \vE_\infty) \to S$ be an admissible $2m$-manifold,
where $S$ is CSC. Then for any admissible K\"ahler class $\Omega$ on $M$, the
extremal polynomial $F_\Omega$ is the unique polynomial $F$ of degree $\leq m
+ 2$ satisfying \eqref{extremal-deriv}--\eqref{extremal-values}
and~\eqref{BCs}.
\end{prop}
\begin{proof} An admissible K\"ahler class on $M$ is specified by
parameters $x_a$ such that $x_0=1$, $x_\infty=-1$ and otherwise $0<|x_a|<1$
with $\omega_a/x_a$ positive. We write $\Scal_{\pm g_a}=\pm 2d_as_a$.
Equation~\eqref{extremal-values} can be solved for a degree $N-1$ polynomial
$P_0$ by Lagrange interpolation, i.e.,
\begin{equation*}
P_0(z)= \sum_{a} 2 d_as_ax_a\prod_{b\in\smash{\hat\cA}\vphantom{I}} (1+x_b z)
\end{equation*}
and then we can write the general degree $N+1$ solution as
\begin{equation*}
P(z) = P_0(z) + (A z + B) \prod_{a\in\smash{\hat\cA}\vphantom{I}} (1+x_a z)
=\biggl(A z+ B+ \sum_{a} \frac{ 2 d_a s_a x_a}{1+x_a z}\biggr)
\prod_{a\in\smash{\hat\cA}\vphantom{I}} (1+x_a z)
\end{equation*}
so that
\begin{equation*}
F''(z) = \biggl(A z+ B+ \sum_{a} \frac{ 2 d_a s_a x_a}{1+x_a z}\biggr)\Mpc(z).
\end{equation*}
Integrating $F''(z)$ and $z F''(z)$ on $[-1,1]$,~\eqref{BCs} now implies that
$A,B$ satisfy~\eqref{system}. Hence $F=F_\Omega$ is the unique solution.
\end{proof}
\begin{rem} An alternative approach is to solve the initial value problem
(at $z=-1$) for $F(z)$. The boundary conditions at $z=1$ then show that $A,B$
satisfy~\eqref{system}. This gives another formula for the extremal
polynomial:
\begin{equation}\label{IVP-formula}
F_\Omega(z) = 2(z+1) \Mpc(-1) + \int_{-1}^z \biggl(A t+ B+ \sum_{a} \frac{
2 d_a s_a x_a}{1+x_a t} \biggr) \Mpc(t) (z-t) dt,
\end{equation}
where $A$ and $B$ are given (as usual) by~\eqref{system}.
\end{rem}

Proposition~\ref{extremal-poly} shows that the existence of an
{\it admissible} extremal K\"ahler metric in $\Omega$ is equivalent to the
positivity of the extremal polynomial $F_\Omega$ on $(-1,1)$. Since the
leading coefficient is a nonzero multiple of $A$, Proposition~\ref{futaki-ext}
shows that such a metric will be CSC iff the Futaki invariant $\mathfrak
F_\Omega(K)$ vanishes.

\begin{rem}\label{openness}
Since $F_\Omega$ depends continuously (in fact analytically) on the admissible
K\"ahler class, it is positive on $(-1,1)$ for an open subset of such classes.
This observation fits in with the general stability result of LeBrun and
Simanca~\cite{Le-Sim2}.
\end{rem}
We now show that $F_\Omega$ is positive on $(-1,1)$ for sufficiently small
$\Omega$.

\begin{prop}\label{adm-exist} Let $M=P(\vE_0\oplus \vE_\infty)\to S$ be
admissible, where $S$ is a local K\"ahler product of CSC metrics $(\pm
g_a,\pm\omega_a)$.  Then there is a nonempty open subset of admissible
K\"ahler classes on $M$ which contain an \textup(admissible\textup) extremal
K\"ahler metric of positive scalar curvature. The admissible K\"ahler classes
containing a CSC metric form a real analytic hypersurface which is nonempty if
$\cb_1(\vE_\infty)-\cb_1(\vE_0)$ is strictly indefinite over $S$ \textup(i.e.,
the definite forms $\omega_a$ do not all have the same sign\textup).
\end{prop}
\begin{proof} 
As we noted in Remark~\ref{openness}, the extremal polynomial $F_\Omega$ is
positive on $(-1,1)$ for an open subset of admissible K\"ahler classes. It
remains to see that this open subset is nonempty and to find the CSC metrics
in the family.  For this, we study the behaviour of $F_\Omega$ near $x_a=0$
for all $a\in \cA$.

\begin{lemma}\label{asymptotics} The coefficients $A$ and $B$ defined by
\eqref{system}, as functions of $x_a\, (a\in \cA)$ for $|x_a|$ small are
given by
\begin{align}\label{A-limit}
A &= -2(2+d_0+d_\infty){\textstyle\sum}_{a\in\cA} d_a x_a +O(x^2)\\
B &= -(1+d_0+d_\infty)(2+d_0+d_\infty)\\
\notag &\qquad -2{\textstyle\sum}_{a\in\cA} d_a s_a x_a
+2(d_0-d_\infty){\textstyle\sum}_{a\in\cA} d_a x_a +O(x^2)
\end{align}
where $O(x^2)$ is shorthand for $\sum_{a,b\in \cA} O(x_ax_b)$.
\end{lemma}
The proof is given in Appendix~\ref{appB}. It follows that in the limit
$x_a\to 0$ for all $a\in \cA$ (which does not give a K\"ahler class),
$F_\Omega''(z)/(1+z)^{d_0} (1-z)^{d_\infty}$ is given by
\begin{equation}\label{inflec}
-(1+d_0+d_\infty)(2+d_0+d_\infty)
+ \frac{2d_0(d_0+1)}{1+z} + \frac{2d_\infty(d_\infty+1)}{1-z}.\\
\end{equation}
If $d_0=0$ and $d_\infty=0$, this is negative on $(-1,1)$ and $F_\Omega$ is
convex. By~\eqref{BCs}, for some $\eps>0$, $F_\Omega$ is positive and
increasing on $(-1,-1+\eps)$ and concave if $d_0>0$, while it is positive and
decreasing on $(1-\eps,1)$ and concave if $d_\infty>\infty$.
By~\eqref{inflec}, for sufficiently small $\Omega$, $F_\Omega$ does not have
enough inflection points to have a zero on $(-1,1)$ and so it is positive
there. Hence the set of admissible K\"ahler classes containing an admissible
extremal K\"ahler metric is nonempty. Since $z\in[-1,1]$, we see that for
$x_a$ sufficiently small, the scalar curvature $- A z - B$ of $g$ is positive.

Now the Futaki invariant $\mathcal F_\Omega(K)\sim A$ is a rational function
of $x_a$, $a\in \cA$, so the CSC metrics form a real analytic hypersurface.
It is then clear from~\eqref{A-limit} that if $x_a$ occur with both signs, $A$
has nonconstant sign for small $x_a$.
\end{proof}
\noindent (Using the Matsushima--Lichnerowicz criterion, the CSC existence
result in this proposition provides an alternative proof of
Proposition~\ref{maximal}.)

\begin{prop}\label{ext=>adm} Let $M$ be an admissible projective bundle over a
CSC base. Then if an admissible K\"ahler class on $M$ contains an extremal
K\"ahler metric, this extremal K\"ahler metric is admissible up to
automorphism.
\end{prop}
\begin{proof} Consider the set $U$ of admissible K\"ahler classes that
contain an extremal K\"ahler metric invariant under the maximal compact
subgroup $G$ of $H_0(M)$ defined in Proposition~\ref{maximal}. By
LeBrun--Simanca~\cite{Le-Sim2}, $U$ is open in the set of all admissible
K\"ahler classes.  Suppose there is some admissible class $\Omega_0$ (with
parameters $x_a^0$) which contains an inadmissible extremal metric $g_0$; by
the Calabi Theorem, we can assume that $g_0$ is $G$-invariant, i.e., $\Omega_0
\in U$.  By LeBrun--Simanca~\cite{Le-Sim2}, this implies that in all K\"ahler
classes sufficiently close to $\Omega_0$ there are $G$-invariant extremal
K\"ahler metrics close to $g_0$ (in suitable Sobolev spaces; by the Sobolev
embedding theorem this also holds in the $C^{\ell}(M)$ topology, for any $\ell
>0$). We have two cases:
\begin{bulletlist}
\item there is an open neighbourhood of $\Omega_0$ in $U$ for which the
extremal polynomial is not positive on $(-1,1)$;
\item there are admissible extremal K\"ahler metrics in classes arbitrarily
close to $\Omega_0$.
\end{bulletlist}
The first case contradicts the existence of admissible extremal K\"ahler
metrics on $M$ for sufficiently small $\Omega$, i.e., the positivity of the
extremal polynomial. Indeed, for such $\Omega$, $F_\Omega(z)$ has at most two
inflection points in $(-1,1)$ by~\eqref{inflec}, and an easy case by case
analysis (according to whether $d_0,d_\infty$ are zero or positive) then shows
that $Q_\Omega(z):=F_\Omega(z)/(1+z)^{\smash{d_0}}(1-z)^{\smash{d_\infty}}$,
as a polynomial in $z$, has simple roots.  However, by assumption, for all
$\Omega$ in some open neighbourhood of $\Omega_0$, the extremal polynomial
$F_\Omega$ is nonnegative (by Corollary~\ref{F-nonneg}) but not positive on
$(-1,1)$, so $Q_\Omega(z)$ has zero discriminant. Since it is analytic in
$\Omega$, it is identically zero, a contradiction.

In the second case we apply instead the uniqueness result~\cite{CT,CT2} for
extremal K\"ahler metrics, as in~\cite{christina1}. Let $\Omega_k$ be a
sequence of admissible K\"ahler classes (with parameters $x_a^k$) which
converges to $\Omega_0$ (i.e., $x_a^k$ converges to $x_a^0$ for all $a$) and
such that $\Omega_k$ contains an admissible extremal K\"ahler metric $\tilde
g_k$ which is not CSC.  By LeBrun--Simanca~\cite{Le-Sim2}, it follows that for
$k\gg 1$ there are $G$-invariant extremal K\"ahler metrics $g_k \in \Omega_k$
which converge to $g_0$ in the $C^2(M)$ topology.  By
Chen--Tian~\cite{CT,CT2}, $g_k$ is the pullback of $\tilde g_k$ by an
automorphism $\Psi_k$ of $(M,J,\Omega_k,G)$. We now claim that $g_0$ is the
pullback by an automorphism of an admissible extremal K\"ahler metric in
$\Omega_0$, completing the proof.

To prove the claim, we use the theory of hamiltonian $2$-forms of order $1$
from~\cite{ACG2,ACGT3}. Since $\tilde g_k$ admits such a $2$-form with $S^1$
action generated by $K$, so does $g_k$ (by Proposition~\ref{futaki-ext}, $K$
is a nonzero multiple of the extremal vector field of $(\Omega_k,G)$ and so is
preserved by $\Psi_k$). Now if $(g,\omega)$ is any K\"ahler metric on $M$ for
which $K=J\grad_g z$ generates an isometric $S^1$ action, then it follows from
\cite{ACGT3} that this action comes from a hamiltonian $2$-form if and only if
it is {\it rigid} (meaning that $g(K,K)$ depends only on $z$) and {\it
semisimple} (meaning that for any regular value $z_0$ of $z$, the
$z$-derivative at $z=z_0$ of the family of K\"ahler quotient metrics $g_{\hat
S}(z)$ on the complex quotient $\smash{\hat S}$ is parallel and diagonalizable
with respect to $g_{\hat S}(z_0)$). Thus, the $S^1$ action generated by $K$ is
rigid and semisimple with respect to $g_k$, hence also with respect to $g_0$
by continuity, so that $g_0$ itself admits a hamiltonian $2$-form of order $1$
with $S^1$ action generated by $K$.

We now apply Theorem~\ref{ACGTthm} and Proposition~\ref{equations} to $g_0$:
it follows that $g_0$ is adapted to the bundle structure of $M=P(E_0\oplus
E_{\infty})\to S$ and induces a CSC K\"ahler metric ${g_{\hat S}}^0$ on
$\smash{\hat S}= P(E_0)\times_{S}P(E_{\infty})$, and a connection $1$-form
$\theta_0$ on the principal $\C\mult$-bundle over $\smash{\hat S}$ (whose
total space is identified with $M^0$), such that the $(1,1)$-form
$d\theta_0=:{\omega_{\hat S}}^0$ is parallel and diagonalizable with respect
to ${g_{\hat S}}^0$. Then (by Chern--Weil theory) we have $[{\omega_{\hat
S}}^0]= [\omega_{\hat S}]$. It follows that the K\"ahler form of ${g_{\hat
S}}^0$ is in the cohomology class $\sum_a [\omega_a]/x_a^0$, since, as
explained in \S\ref{admissible-class}, this cohomology class is determined by
the admissible class $\Omega_0$, which can be uniquely written as the sum of
the `projective Thom class' $\smash{\hat\Xi}$, of $\smash{\hat M}= P(\cO
\oplus \hat L)\to \hat S$ and a pullback from $\hat S$.

Now by Chen--Tian~\cite{CT,CT2} again, there is an automorphism $\psi$ of
$\smash{\hat S}$ with $g_{\hat S}=\psi^*{g_{\hat S}}^0$, since these are CSC
K\"ahler metrics in the same K\"ahler class. In fact, the proof of
\cite{CT,CT2} essentially shows that any two extremal metrics in a given
K\"ahler class can be connected by a geodesic in the space of K\"ahler
potentials, and therefore $\psi$ can be chosen in the reduced (connected)
automorphism group $H_0(\smash{\hat S})\subseteq {\Aut}_0(\smash{\hat S})$
(see e.g.~\cite{guan1, gauduchon}); in particular, such a $\psi$ acts
trivially on cohomology. By Proposition~\ref{isometry algebra} there is a
fibre-preserving $S^1$-equivariant automorphism $\Psi$ of $\smash{\hat M}=
P(\cO \oplus \hat L)\to \smash{\hat S}$, which induces $\psi$ on $\smash{\hat
S}$. Thus, $\Psi$ preserves the $\C\mult$-bundle structure of $M^0 \to \hat S$
and sends the connection $1$-form $\theta_0$ to a connection $1$-form $\tilde
\theta_0$ with curvature $d{\tilde \theta}_0= \psi^*{\omega_{\hat S}}^0\in
[{\omega_{\hat S}}^0]= [\omega_{\hat S}]$; now since $\psi^*{\omega_{\hat
S}}^0$ and $\omega_{\hat S}$ are both parallel (and therefore harmonic) with
respect to $\psi^* {g_{\hat S}}^0 = g_{\hat S}$, Hodge theory implies they are
equal.  We can therefore send $\tilde \theta_0$ to $\theta$ via a bundle
isomorphism.

Thus we have constructed (on $M^0$ and hence, by a standard extension
argument, everywhere) an automorphism sending $g_0$ to an admissible extremal
K\"ahler metric in $\Omega_0$, as required.
\end{proof}
Theorem~\ref{main-thm} follows from
Propositions~\ref{extremal-poly}--\ref{ext=>adm}.

\section{Existence and nonexistence results for extremal K\"ahler metrics}
\label{existence?}

In this section we use Theorem~\ref{main-thm} to construct explicit examples
of extremal K\"ahler metrics. We also obtain some nonexistence results for CSC
K\"ahler metrics.

\subsection{Constructing admissible extremal K\"ahler metrics}
\label{general-theory}

We begin with a root counting argument due to Hwang~\cite{hwang} and
Guan~\cite{guan0} which gives a complete construction when the base $S$ is a
local K\"ahler product of {\it nonnegative} CSC K\"ahler metrics (in fact
Hwang and Guan only considered the case that $S$ has constant nonnegative
eigenvalues of the Ricci tensor, but the proof is no different in general, and
the idea to weaken this hypothesis is already explored
in~\cite{hwang-singer}).

\begin{prop}\label{root-count}
Suppose that $M=P(E_0\oplus E_\infty)\to S$ is admissible where $S$ is a local
K\"ahler product of nonnegative CSC metrics. Then every admissible K\"ahler
class contains an \textup(admissible\textup) extremal K\"ahler metric.
\end{prop}
\begin{proof}
By the boundary conditions $F_\Omega$ is positive, and increasing or
decreasing, on $(-1,-1+\eps)$ or $(1-\eps,1)$ respectively, for some $\eps>0$.
Suppose it is not positive on $(-1,1)$. Then it has at least two maxima, one
minimum and two inflection points on $(-1,1)$. It follows that $P$ has at
least two roots in $(-1,1)$.

Let $y_1\leq \cdots \leq y_Q$ and $z_1\leq \cdots \leq z_R$ ($Q,R\geq 0$)
denote the roots (counted with multiplicity) of $P$ in $[1,\infty)$ and
$(-\infty,-1]$ respectively, and put $y_0=1$, $y_{Q+1}=\infty$, $z_0=-\infty$,
$z_{R+1}=-1$. We order $\{x_a:a\in\smash{\hat\cA}\}$ as
\begin{equation*}
-1\leq x_{a_1}<\cdots < x_{a_J} < 0 < x_{a_{J+1}} < \cdots < x_{a_N}\leq 1
\end{equation*}
(for some $0\leq J\leq N$) so that $g_{a_j}$ is negative definite (hence with
$s_{a_j}$ nonpositive) for $j\leq J$ and positive definite (hence with
$s_{a_j}$ nonnegative) for $j\geq J+1$.

Therefore by~\eqref{extremal-values}, for each $0\leq q\leq Q$, there is at
most one $x_{a_j}$ with $y_q\leq -1/x_{a_j}< y_{q+1}$, so that $Q+1\geq J$
with equality iff there is exactly one $x_{a_j}$ in each such interval.
Similarly, for each $0\leq r\leq R$, there is at most one $x_{a_j}$ with $z_r<
-1/x_{a_j}\leq z_{r+1}$, so that $R+1\geq N-J$ with equality iff there is
exactly one $x_{a_j}$ in each such interval. Thus $P$ has at least $N-2$ roots
outside $(-1,1)$.

Since $P$ has degree $\leq N+1$, it has at most $N-1$ roots outside $(-1,1)$,
so we must either have $Q+1=J$ or $R+1=N-J$. If (without loss of generality)
$Q+1=J$ then $-1/x_{a_1}<y_1$, so that $P(-1/x_{a_1})>0$
(by~\eqref{extremal-values} again) and there must be a root of $P$ between
$-1/x_{a_1}$ ($1 \leq -1/x_{a_1} < y_1 $) and the last maximum of
$F_\Omega$ in $(-1,1)$. This now forces $R+1=N-J$ also, hence
$P(-1/x_{a_N})>0$ and there must be a root of $P$ between the first maximum of
$F_\Omega$ in $(-1,1)$ and $-1/x_{a_N}$ ($z_R < -1/x_{a_N} \leq -1 $),
contradicting $\deg P\leq N+1$.
\end{proof}

Because of this result, in the rest of this section we shall mainly be
interested in the influence of negative scalar curvature factors in the base
metrics. In the presence of such factors, the existence of extremal K\"ahler
metrics is nontrivial, as was already observed in~\cite{christina1} for ruled
surfaces. By Theorem~\ref{main-thm} such a metric exists in a given admissible
class $\Omega$ iff the extremal polynomial $F_\Omega$ is positive on $(-1,1)$.
However, the integrals $\alpha_i$ and $\beta_j$ involved in the above
construction of $F_\Omega$ are hard to compute in general (see
Appendix~\ref{appB}).  For the next examples, we therefore adopt a different
approach to compute $F_\Omega$.  Instead of solving
\eqref{extremal-deriv}--\eqref{extremal-values} and integrating, we solve
first the boundary conditions.

It is easy to see that~\eqref{BCs} together with~\eqref{extremal-values} for
$a=0$ and $a=\infty$ (if there are blow-downs) are solved by any $F$ of the
form
\begin{equation}
F(z) = (1-z^2)\bigl(\Mpc(z)+(1+z)^{d_0+1}(1-z)^{d_\infty+1} q(z)\bigr)
\end{equation}
for some polynomial $q(z)$. Conversely any polynomial solution is of this
form, and to obtain the extremal polynomial $F_\Omega$, the degree of $q$
must be $\leq (\sum_{a\in\cA} d_a)-1$. Now it remains to compute $F''(z)$, to
solve~\eqref{extremal-deriv}--\eqref{extremal-values} for $a\in\cA$ (so that
$F=F_\Omega$), and to check positivity.  For a given projective bundle and
admissible K\"ahler class this leads to equations on the coefficients of
$q$. The Futaki invariant will be zero (and the metric will be CSC) iff $q$
has degree $\leq (\sum_{a\in\cA} d_a)-2$.

In general, the algebraic equations on $q$ are hopelessly complicated.
However, when $S$ has real dimension $\leq 4$, they are tractable.

\subsection{Extremal K\"ahler metrics over a Riemann surface}
\label{extRS}

We consider first extremal K\"ahler metrics on projective bundles over a
Riemann surface, generalizing the study of projective line bundles
in~\cite{christina1}.

Let $\Sigma$ be a compact Riemann surface with CSC metric $(\pm
g_\Sigma,\pm\omega_\Sigma)$ and let $M = P(\vE_0\oplus \vE_\infty) \to
\Sigma$, where $\vE_0$, $\vE_\infty$ are projectively-flat hermitian vector
bundles with ranks $d_0+1>0$, $d_\infty + 1> 0$, and
$\cb_1(\vE_\infty)-\cb_1(\vE_0) =[\omega_\Sigma/2\pi]$.  Let $\pm 2s$ be the
scalar curvature of $\pm g_{\Sigma}$ and $\Omega$ be an admissible K\"ahler
class on $M$ defined by $0<|x|<1$.  The metric
\begin{equation*}
g = (z+1) g_0 + (z + 1/x) g_\Sigma +(z-1)g_\infty + \frac{\Mpc(z)}{F(z)}
dz^2 + \frac{F(z)}{\Mpc(z)} \theta^2,
\end{equation*}
is an extremal metric in the given class iff $F=F_\Omega$. From
\S\ref{general-theory}, we know that
\begin{equation*}
F_\Omega(z)=(1+z)^{d_0+1}(1-z)^{d_\infty+1}\bigl((1+xz) + c(1-z^2)\bigr)
\end{equation*}
where $q(z)=c$ is a constant uniquely determined by the equation
\begin{equation*}
F_\Omega''(-1/x) = 2sx(1-1/x)^{d_0}(1+1/x)^{d_\infty}
\end{equation*}
(this holds whether or not $d_0,d_\infty$ are zero).  We solve this to obtain
\begin{equation*}
c(s,x)= -\tfrac{2x^2 (2+d_0(1+x)+d_\infty(1-x)-sx)}
{(2+d_0(1+x)+d_\infty(1-x))(4+d_0(1+x)+d_\infty(1-x))
+(4+d_0+d_\infty)(1-x^2)}.
\end{equation*}
Since $sx$ has same sign as $\Scal_{\Sigma}$, it can be positive only when
$\Sigma= \C P^1$, in which case $E_0 = \cL_0 \otimes \C^{d_0+1}$ and
$E_{\infty} = \cL_{\infty}\otimes \C^{d_{\infty}+1}$ for some line bundles
$\cL_{0}, \cL_{\infty}$. It then follows that $\omega_\Sigma$ is integral, and
thus $s=p/q$ where $p \leq 2$ and $q$ is an integer of same sign as $x \in
(-1,1)\setminus \{0 \}$ (see Remark \ref{integrality}); then we have that
$sx<2$, so that $c<0$.  Therefore $\mathfrak F_\Omega(K)$, which is a
nonzero multiple of $c$, doesn't vanish for any admissible K\"ahler class.
Since $b_2(\Sigma)=1$, every K\"ahler class on $M$ is admissible, so we get an
immediate nonexistence result.

\begin{thm}\label{noCSCsigma}
Let $\vE_0$, $\vE_\infty$ be projectively-flat hermitian vector bundles over a
Riemann surface $\Sigma$. Then there are no CSC K\"ahler metrics on $M=
P(\vE_0\oplus \vE_\infty)$ unless $\cb_1(\vE_0)=\cb_1(\vE_\infty)$
\textup(i.e., $\vE_0\oplus \vE_\infty$ is polystable\textup).
\end{thm}
This partially extends the converse in Corollary~\ref{rscase} to the case
${\mathbf g}>1$. Compared to Theorem~\ref{main-thm}, we note that here
$\cb_1(\vE_\infty)-\cb_1(\vE_0)$ can never be strictly indefinite. On the
other hand, by Theorem~\ref{main-thm}, we have an extremal K\"ahler metric for
sufficiently small $\Omega$. Indeed it is easy to see that $|c|$ is small when
$|x|$ is small, and hence $(1+xz) + c(1-z^2)$ is positive on $(-1,1)$. We also
know from Proposition~\ref{root-count} that if $\Sigma$ has genus $0$ or $1$,
then every admissible K\"ahler class contains an extremal K\"ahler metric. Let
us now see what happens when ${\mathbf g} > 1$, i.e., when $sx < 0$.

Since $c<0$, the quadratic $Q(z)=(1+xz) + c(1-z^2)$ is concave. It is clearly
positive at $z=\pm 1$, so it is positive on $(-1,1)$ unless its minimum is in
$(-1,1)$ and it is nonpositive there. The minimum value $1+c+x^2/4c$ occurs at
$z=x/2c$ and
\begin{equation*}
c(s,1)=  -\frac{2(1+d_0)-s}
{2(1+d_0)(2+d_0)},\qquad c(s,-1)=  -\frac{2(1+d_\infty)+s}
{2(1+d_\infty)(2+d_\infty)}.
\end{equation*}
It follows that if $s<-d_0(d_0+1)$, then $c(s,1)<-\frac12$ and hence for
$0<x<1$ sufficiently close to $1$, we have $c(s,x)<-\frac12$ and the minimum
of $Q(z)$ is in $(-1,0)$ and nonpositive.  Similarly, if
$s>d_\infty(d_\infty+1)$ then for $-1<x<0$ sufficiently close to $-1$, we have
$c(s,x)<-\frac12$ and the minimum of $Q(z)$ is in $(0,1)$ and nonpositive.

Hence if $\Sigma$ has genus ${\mathbf g}>1$ and $s<-d_0(d_0+1)$ or
$s>d_\infty(d_\infty+1)$ then not every admissible K\"ahler class contains an
admissible extremal metric.

Conversely if $s\geq -d_0(d_0+1)$ then it is easy to check that $c(s,x)\geq
-x/2$ for all $0<x<1$, so the minimum of $Q(z)$ is not in $(-1,1)$ for any
such $x$, whereas if $s\leq d_\infty(d_\infty+1)$, $c(s,x)\geq x/2$ for all
$-1<x<0$ and again the minimum of $Q(z)$ is not in $(-1,1)$ for any such $x$.

\begin{thm}\label{CSCsigma}
Let $\vE_0$, $\vE_\infty$ be projectively-flat hermitian vector bundles of
ranks $d_0+1$, $d_\infty+1$ over a compact Riemann surface $\Sigma$ of genus
${\mathbf g}$, and suppose $\cb_1(\vE_\infty)-\cb_1(\vE_0)=
[\omega_\Sigma/2\pi]$ for a K\"ahler form $\pm \omega_\Sigma$ of constant
curvature. Then there exist admissible extremal K\"ahler metrics on
$P(\vE_0\oplus\vE_\infty)\to \Sigma$.  Such metrics exist in every K\"ahler
class if ${\mathbf g} =0$ or $1$.  For ${\mathbf g} >1$, put
$\rho_\Sigma=s\omega_\Sigma$. Then such metrics exist in every K\"ahler class
if and only if $-d_0(d_0+1)\leq s\leq d_\infty(d_\infty+1)$, otherwise such
metrics exist for $|x|$ sufficiently small \textup(depending on $s$\textup).
\end{thm}

\begin{rem}\label{general-christina} In absence of blow-downs, we 
recover the examples of~\cite{christina1} on (complex) pseudo-Hirzebruch
surfaces $P(\cO\oplus\cL)\to\Sigma$, where there are K\"ahler classes which do
not contain an extremal K\"ahler metric (if $\Sigma$ has genus $\mathbf
g>1$). Our result extends these examples to higher rank projective bundles.
However, in the presence of blow-downs, there do exist projective bundles for
which there is an extremal K\"ahler metric in every K\"ahler class, even with
$\mathbf g>1$.
\end{rem}

\subsection{Nonexistence of CSC K\"ahler metrics over a Hodge $4$-manifold}
\label{extHodge}

We now obtain a similar nonexistence result to Theorem~\ref{noCSCsigma} when
$\dim S=4$.
\begin{thm}\label{noCSCsurface}
Let $(S,\pm g_S,\pm \omega_S)$ be a CSC Hodge $4$-manifold and let $\vE_0$,
$\vE_\infty$ be projectively-flat hermitian vector bundles of ranks $d_0+1$,
$d_\infty+1$ over $S$ with $\cb_1(\vE_\infty)-\cb_1(\vE_0) = [\omega_S/2\pi]$.
Then there are no CSC K\"ahler metrics in the admissible K\"ahler classes on
$M=P(\vE_0\oplus\vE_\infty) \to S$.
\end{thm}
\begin{proof}
We will prove that $\mathfrak{F}_{\Omega}(K)$ is nonzero for any admissible
K\"ahler class $\Omega$ by showing that the leading coefficient of
$F_\Omega$ cannot vanish. Following the discussion in \S\ref{general-theory}
we see that for given $d_0,d_{\infty}\geq 0$ and $0 < |x| < 1$ we have
\begin{equation*}
F_\Omega(z) = (1+z)^{d_0+1}(1-z)^{d_\infty+1}\bigl((1+x z)^{2} + (c z 
+ e)(1-z^2)\bigr)
\end{equation*}
with $c$ and $e$ being constants uniquely determined by the conditions
\begin{equation*}
F_\Omega''(-1/x) = 0 \qquad\text{and} \qquad
\frac{F_\Omega''(z)}{1+xz}\Big|_{z=-1/x} = 4 s x 
(1-1/x)^{d_{0}}(1+1/x)^{d_{\infty}}.
\end{equation*}
The leading coefficient of $F_\Omega$ vanishes iff $c=0$. As before $c$
(and $e$) are determined by $s$ and $x$. In particular, $c(s,x)
={n(s,x)/d(s,x)}$, where
\begin{align*}
d(s,x) &= (1+d_\infty) (2+d_\infty)^2 (3+d_\infty) (1-x)^4\\
&\,+ 4 (2+d_0) (1+d_\infty) (2+d_\infty) (3 + d_\infty) (1-x)^3 (1+x)\\
&\,+ 6 (2+d_0) (2+d_\infty) \bigl(4+d_0+d_\infty+(1+d_0) (1+d_\infty)\bigr)
(1-x)^2(1+x)^2\\
&\,+ 4 (1+d_0) (2+d_0) (3+d_0) (2+d_\infty) (1 - x) (1 + x)^3\\
&\,+ (1+d_0) (2+d_0)^2 (3+d_0) (1 + x)^4
\end{align*}
(which is manifestly positive for $|x|<1$) and
\begin{align*}
-n(s,x)/2x^{3} &=  4(6-3sx+sx^{3})\\
&\, + d_{0}(1+x)\bigl((5-x)(1+x)+(7-x)(3-sx)\bigr)\\
&\, + d_{\infty}(1-x)\bigl((5+x)(1-x)+(7+x)(3-sx)\bigr)\\
&\, +(9-sx) \bigl(d_{0}(1+x)  + d_{\infty}(1-x)\bigr)^{2}
+ \bigl(d_{0}(1+x)+ d_{\infty}(1-x)\bigr)^{3}.
\end{align*}
Since $s=p/q$ where $p\leq 3$ (see Remark~\ref{integrality}) and
$q$ is an integer of same sign as $x \in (-1,1)\setminus \{0 \}$, we have
that $sx < 3$ and a moment's thought then gives that $n(s,x)$, and therefore
$c(s,x)$, is never zero.
\end{proof}

\subsection{CSC K\"ahler metrics over a product of two Riemann surfaces}

As counterpoint to the nonexistence results of
\S\S\ref{extRS}--\ref{extHodge}, we now explore explicitly the existence of
CSC K\"ahler metrics, given by Theorem~\ref{main-thm}, in the simplest case
when the base is a global product of two Riemann surfaces and there are no
blow-downs.

Let $\Sigma_a$ $(a=1,2)$ be compact Riemann surface with CSC metrics $(\pm
g_a,\pm \omega_a)$ and let $M$ be $P(\cO\oplus\cL)\to \Sigma_1\times\Sigma_2$
where $\cL=\cL_1\otimes\cL_2$ and $\cL_a$ are pullbacks of line bundles on
$\Sigma_a$ with $c_1(\cL_a)=[\omega_a/2\pi]$. Let $\pm 2s_a$ be the scalar
curvature of $\pm g_a$ and $-1/x_a$ be the constant roots defining an
admissible K\"ahler class with $x_1 \neq x_2$ (the case $x_1=x_2$ was
considered in \S\ref{extHodge}, where we established nonexistence of CSC
metrics).  We thus have $\Mpc(z)=(1+x_1 z)(1+x_2 z)$ and the metric becomes
\begin{equation*}
g = \frac{1+x_1z}{x_1} g_1 +\frac{1+x_2z}{x_2}g_2 +
\frac{\Mpc(z)}{F(z)} dz^2 + \frac{F(z)}{\Mpc(z)} \theta^2.
\end{equation*}
According to \S\ref{general-theory}, to obtain a CSC metric, $F(z)$ must be
the extremal polynomial
\begin{equation*}
F_\Omega(z)=(1-z^2)\bigl((1+x_1z)(1+x_2z) + cx_1x_2(1-z^2)\bigr)
\end{equation*}
where $c$ is a constant such that the following relations are satisfied
\begin{equation*}
F_\Omega''(-1/x_1) = 2s_1(x_1-x_2),\qquad
F_\Omega''(-1/x_2) = 2s_2(x_2-x_1).
\end{equation*}
Writing $2(1-c)=s$ (which is $\frac16\Scal_g$ and not to be confused with the
$s$ in the previous two paragraphs), these relations hold iff
\begin{align}\label{CSCeq1}
x_1(s_1(x_1-x_2)-2+(1-s)x_1x_2)+3(s-1)x_2&=0\\
x_2(s_2(x_2-x_1)-2+(1-s)x_1x_2)+3(s-1)x_1&=0,
\label{CSCeq2}
\end{align}
and these are precisely the conditions on an admissible K\"ahler class
$\Omega$ (parameterized by $(x_1,x_2)$ with $0<|x_a|<1$) coming from the
vanishing of $\mathfrak{F}_{\Omega}(K)$ (see
Proposition~\ref{extremal-poly}).  Eliminating $s=2(1-c)$, we obtain (using
$x_1 \neq x_2$)
\begin{equation}\label{hypersurface}
x_{1}(6 + s_{1}x_{1}(x_{2}^{2}-3)) + x_{2}(6+s_{2}x_{2}(x_{1}^{2}-3)) =0.
\end{equation}

The normalized scalar curvatures $s_{a}$ are subject to the integrality
conditions $s_{a} =2(1-{\mathbf g}_a)/q_a$ for $q_a$ a nonzero integer with
the same sign as $x_a$, where ${\mathbf g}_a$ is the genus of $\Sigma_{a}$. In
particular $s_a x_a<2$.  This and \eqref{hypersurface} imply that $x_{1}x_{2}
< 0$; we thus get a nonexistence result in the case $x_1x_2>0$.

\begin{thm}\label{no-CSC-product-surfaces} Let $(\Sigma_a, \omega_a)$
\textup($a=1,2$\textup) be compact CSC Riemann surfaces and $\cL$ be a
holomorphic vector bundle over $\Sigma_1\times \Sigma_2$ with $c_1(\cL) =
[(\omega_1+\omega_2)/2\pi]$ \textup(so that $c_1(\cL)$ is positive
definite\textup). Then there are no admissible K\"ahler classes on $M= P(\cO
\oplus \cL) \cong P(\cO\oplus \cL^{-1})$ containing a CSC K\"ahler metric.
\end{thm}

\begin{rem}\label{rem4}
Note that we do not need to assume that the base is a global product of
compact Riemann surfaces for the nonexistence result in the above theorem. It
is sufficient to have a compact base $S$ that is a local product of Riemann
surfaces with CSC and $s_{a}x_{a}<2$, which is always satisfied, since
$\Scal_{\pm g_a}\leq 4$ by the integrality of the pull-back of $\pm\omega_a$
to the universal cover of $S$.
\end{rem}

In contrast to this result, we have the following observation.

\begin{lemma}\label{CSC-product-surfaces} Let $\Omega$ be an admissible
K\"ahler class, corresponding to a solution $(x_1,x_2)$ of
\eqref{CSCeq1}--\eqref{CSCeq2} with $s\geq 0$. Then $\Omega$ admits an
admissible CSC K\"ahler metric with scalar curvature $6s$.
\end{lemma}
\begin{proof}
If \eqref{hypersurface} holds, the extremal polynomial $F_\Omega$ of an
admissible K\"ahler class gives rise to a globally defined CSC K\"ahler metric
iff $F_\Omega>0$ on $(-1,1)$.  Let $Q(z)= F_\Omega(z)/(1-z^2)$, and
observe that the coefficient of $z^2$ in this quadratic is $\frac12
sx_1x_2$. Since $Q(\pm 1)>0$, $Q$ will be positive on $[-1,1]$ if it is
convex, i.e., if $sx_1x_2<0$. If $s=0$, $Q(z)$ is linear and positive on
$[-1,1]$. Since $x_{1}x_{2}<0$, $F_\Omega$ is positive on $(-1,1)$
whenever we have solutions of \eqref{CSCeq1}--\eqref{CSCeq2} with $s\geq 0$.
\end{proof}

We now obtain some explicit solutions of \eqref{CSCeq1}--\eqref{CSCeq2}.  If
we take $x_2\eqref{CSCeq1}-x_1\eqref{CSCeq2}$ and
$x_1\eqref{CSCeq1}-x_2\eqref{CSCeq2}$, we obtain, for $x_1\neq x_2$:
\begin{align*}
(s_1+s_2)x_1x_2&=3(s-1)(x_1+x_2)\\
2(x_1+x_2)&=s_1x_1^2+s_2x_2^2+(1-s)x_1x_2(x_1+x_2).
\end{align*}
These are equivalent to~\eqref{CSCeq1}--\eqref{CSCeq2} for $x_1^2\neq
x_2^2$. As $x_1\neq x_2$, $x_1^2=x_2^2$ iff $x_1+x_2=0$ and then
$s_1+s_2=0$. The following lemma deals with this case.

\begin{lemma} If $s_1+s_2=0$, then either $x_1+x_2=0$ and
$s=(1-x_1^2+2s_1x_1)/(3-x_1^2)$, or, without loss, $x_1=x_2+1$, $s=1$, and
$s_1=2=-s_2$.  Conversely, these give solutions
of~\eqref{CSCeq1}--\eqref{CSCeq2}.
\end{lemma}
\begin{proof} Clearly $s_1+s_2=0$ iff $s=1$ or 
$x_1+x_2=0$. The formula for $s$ in the latter case is immediate
from~\eqref{CSCeq1}. Now if $s=1$, then without loss of generality $s_1=-s_2$
is nonnegative and we must have either $x_1+x_2=0$, or $s_1>0$ and
$x_1=x_2+2/s_1$. Since $0<|x_a|<1$, this forces $x_1$ to be positive, hence
$s_1\leq 2$, so in fact we must have $s_1=2$ and $x_1=x_2+1$.
\end{proof}

In order to apply Lemma~\ref{CSC-product-surfaces}, we suppose in the first
case above that $x_1s_1\geq 0$: then $s>0$ since $1-x_1^2>0$ for $|x_1|<1$.
Thus in both cases $s_ax_a\geq 0$ for $a=1,2$ and we obtain CSC K\"ahler
metrics on projective line bundles over $T^2\times T^2$, $T^2\times \C P^1$
and $\C P^1\times \C P^1$. In particular any K\"ahler class on $P(\cO\oplus
\cO(q,-q))\to \C P^1\times \C P^1$ is admissible, so the above Lemmas and
Proposition~\ref{root-count} yield the following conclusions.
\begin{thm}\label{CSCKoiso-Sakane}
On $P(\cO\oplus \cO(q,-q)) \to \C P^1 \times \C P^{1}$ $(q\geq 1)$,
any K\"ahler class \textup(parameterized, up to scale, by $0<x_1<1$
and $-1<x_2<0$\textup) contains a unique admissible extremal K\"ahler
metric. For $q>1$ this metric is CSC if and only if $x_1+x_2=0$, while
for $q=1$ it is CSC if and only if $x_1+x_2=0$ or $x_1= x_2+1$.
\end{thm}
When $q=1$, the two $1$-parameter families of CSC K\"ahler classes of this
theorem intersect at $x_1=1/2$, $x_2=-1/2$.  In fact, the CSC metric in this
K\"ahler class is the Koiso--Sakane K\"ahler--Einstein metric~\cite{koi-sak1}.

We end our study of CSC K\"ahler metrics on $P(\cO\oplus\cL) \to
\Sigma_1\times \Sigma_2$, by considering the case of zero scalar curvature
metrics, which we do not obtain automatically from Theorem~\ref{main-thm}.
If $s=0$ then equation~\eqref{CSCeq1} defines $x_{2}$ as a function of $x_{1}$
\begin{equation*}
x_2=f_{1}(x_{1})=x_{1}\frac{2-s_{1}x_{1}}{x_{1}^{2}-s_{1}x_{1}-3},
\end{equation*}
whereas~\eqref{CSCeq2} defines $x_{1}$ as a function of $x_{2}$
\begin{equation*}
x_1=f_{2}(x_{2})=x_{2}\frac{2-s_{2}x_{2}}{x_{2}^{2}-s_{2}x_{2}-3}.
\end{equation*}
Note that $f_{1}(0)=f_{2}(0)=0$ and the gradients $dx_{2}/dx_{1}$ of the two
graphs at $x_{1}=x_{2}=0$ are both negative. By comparing the size of the
gradients one sees that for $x_{1}$ small and positive the graph of $f_{1}$ is
above the graph of $f_{2}$. Note also that the denominator appearing in
$f_{a}(x_{a})$ is negative at $x_{a}=0$.

Assume that $s_{1} \leq 0$. If $f_{1}$ has no asymptotes for $0<x_{1}<1$ then
$f_{1}(1) \leq -1$. Otherwise, for the asymptote $x_{1}=v$ closest to
$x_{1}=0$ we have $\lim_{x_{1} \rightarrow v^{-}} = -\infty$.  Assume moreover
that $0 < s_{2}$. If $f_{2}$ has no asymptotes for $-1 < x_{2} <0$ then
$f_{2}(-1) >1$. Otherwise, for the asymptote $x_{2}=v$ closest to $x_{2}=0$ we
have $\lim_{x_{2} \rightarrow v^{+}} = +\infty$.  By continuity the graphs of
$f_{1}$ and $f_{2}$ intersect in the open square $(0,1) \times (0,-1)$ and
\eqref{CSCeq1}--\eqref{CSCeq2} is solved for some $0< x_{1} <1$ and $-1 <
x_{2} <0$.

\begin{thm}
Let $(\Sigma_a, \pm \omega_a)\; (a=1,2)$ be compact Riemann surfaces with
genus $\mathbf g_{a}$ and canonical bundles $\cK_a$, and suppose that the
K\"ahler forms $\pm \omega_a$ are integral with constant curvature.  Let
$\cL_a$ be line bundles on $\Sigma_a$ with $c_1(\cL_a)=[\omega_a/2\pi]$ and,
if $\mathbf g_a \neq 1$, let $\cL_a$ be $\cK_a^{\smash{q_a/2(\mathbf
g_{a}-1)}}$ tensored by a flat line bundle, for an integer $q_a$.

Then, there is an admissible scalar-flat K\"ahler metric on
$P(\cO\oplus\cL_1\otimes\cL_2)\to \Sigma_1\times \Sigma_2$ in the following
cases\textup:
\begin{itemize}
\item $\Sigma_{1}=T^{2}$ and $\cL_1$ is ample, $\Sigma_2$ has genus $\mathbf
g_{2}>1$ and $q_{2} < 0$\textup;
\item $\Sigma_1$ and $\Sigma_2$ both have genus $\mathbf g_{a}>1$,
$q_{1} > 0$, and $q_{2} < 0$.
\end{itemize}
\end{thm}

\section{K-stability and admissible extremal K\"ahler metrics} \label{stab}

\subsection{Introduction to stability}\label{intro}

For a Hodge manifold $(M,\Omega$), the existence of a CSC K\"ahler metric in
$\Omega$ is conjectured to be equivalent to a notion of
{\it stability}~\cite{CT,Do1,Do2,luo,PT,tian,zhang} for the polarized
projective variety $(M,L)$, where $L$ is a line bundle on $M$ with $c_1(L)=
\Omega/2\pi$. This conjecture is drawn from a detailed formal picture which
makes clear an analogy with the well-established relation between the
polystability of vector bundles and the existence of Einstein--Hermitian
connections.

At present the most promising candidate for the conjectured stability
criterion is `K-polystability', in the form given by Donaldson~\cite{Do2},
following Tian~\cite{tian}: a polarized projective variety $(M,L)$ is
K-polystable if any `test configuration' for $(M,L)$ has nonpositive Futaki
invariant with equality iff the test configuration is a product. We shall
explain this definition shortly. We also discuss an idea of Ross and
Thomas~\cite{RT,RT2}, who focus on test configurations arising as
`deformations to the normal cone' of subschemes of $(M,L)$, leading to a
notion of `slope' K-polystability analogous to the slope polystability of
vector bundles. We explore this analogy further in \S\ref{stab-CSC}.

(Note that some authors use the term {\it K-stable} rather than K-polystable,
but the latter term agrees better with pre-existing notions of stability.)

\subsubsection{Finite dimensional motivation}\label{finitedim}

Let $(X,\cL,\Omega)$ be a polarized K\"ahler manifold with a hermitian metric
on $\cL$ with curvature $-i\Omega$ (thus $c_1(\cL)=\Omega/2\pi$).  Suppose a
compact connected group $G$ acts holomorphically on $X$ with momentum map
$\mu\colon X\to \mathfrak g^*$ (i.e.,
$d\langle\mu,\xi\rangle=-\Omega(K_\xi,\cdot)$, where $K_\xi$ is the vector
field on $X$ corresponding to $\xi\in\mathfrak g$, the Lie algebra of
$G$). There is a lift of the action to $\cL$ generated by $\tilde K_\xi+
\langle\mu,\xi\rangle K$ for each $\xi\in \mathfrak g$, where
$\langle\mu,\xi\rangle$ is pulled back to $\cL$, $\tilde K_\xi$ is the
horizontal lift, and $K$ generates the standard $U(1)$ action on $\cL$.  The
action of $\mathfrak g$ on $X$ and $\cL$ extends to an action of the
complexification $\mathfrak g^c$ and we assume this integrates to an action of
a complex Lie group $G^c$.

By a well-known result of Kempf--Ness and Kirwan, for any $x\in X$, there is a
$g\in G^c$ such that $\mu(g\cdot x)=0$ iff for any nonzero lift $\tilde x$ of
$x$ to $\cL^*$, the orbit $G^c\cdot \tilde x=0$ is closed.  Such points $x$
are said to be {\it polystable}. If $X^{ps}$ denotes the set of polystable
points in $X$, we then have an equality between $X^{ps}/G^c$, the polystable
quotient of $X$ by $G^c$, and the symplectic quotient $X/\!/G=\mu^{-1}(0)/G$.

$G^c\cdot\tilde x$ is closed iff $\alpha(\C\mult)\cdot \tilde x$ is closed for
any one parameter subgroup $\alpha\colon\C\mult\mapsinto G^c$. This leads to
the Hilbert--Mumford criterion for polystability: $x$ is said to be
{\it semistable} if for any one parameter subgroup
$\alpha\colon\C\mult\mapsinto G^c$, the linear action of $\C\mult$ on
$\cL_{x_0}^*$ has nonpositive weight $w_{x_0}(\alpha)\leq 0$, where
$x_0=\lim_{\lambda\to 0} \alpha(\lambda)\cdot x$ is the limit point; $x$ is
then polystable if it is semistable and $w_{x_0}(\alpha)=0$ only when $x_0=x$;
finally $x$ is {\it stable} if it is polystable and has zero dimensional
isotropy subgroup.

\subsubsection{The infinite dimensional analogue}\label{infinitedim}

We apply the finite dimensional picture above formally to an infinite
dimensional setting in which $X$ is the space of compatible complex structures
on a compact symplectic manifold $(M,\omega)$ with $H^1(M)=0$. The space $X$
has a natural K\"ahler metric with respect to which the group $G$ of
symplectomorphisms of $M$ acts holomorphically with a momentum map $\mu\colon
X\to C^\infty_0(M,\R)$ given by the scalar curvature of the corresponding
K\"ahler metric on $M$, modified by a constant in order to lie in $\mathfrak
g^*\cong\mathfrak g=C^\infty_0(M,\R)$, the functions with total integral zero,
which is the Lie algebra of the symplectomorphism group equipped with the
$L_2$-inner product.  A quick way to see this is to observe that the Mabuchi
K-energy (see~\S\ref{K-energy-ext}) of $M$ is a K\"ahler potential for the
metric on $X$: the gradient on $X$ of the Mabuchi K-energy is the scalar
curvature~\cite{gauduchon}.

There is no group whose Lie algebra is the complexification $\mathfrak g^c$,
but one can still consider the foliation of $X$ given by the vector fields
induced by $\mathfrak g^c$.  The complex structures in a given leaf are all
biholomorphic by a diffeomorphism in the connected component of the identity,
and pulling back the symplectic form $\omega$ by these biholomorphisms, we may
identify the leaf with the set of all K\"ahler metrics in a fixed K\"ahler
class, compatible with a fixed complex structure on $M$. Hence there should be
a CSC metric in a given K\"ahler class iff the momentum map $\mu$ vanishes on
the corresponding leaf iff the leaf is stable in a suitable sense.

To make precise this infinite dimensional analogue, we formalize what is meant
by the orbit of a $1$-parameter subgroup in terms of `test configurations' and
give a Hilbert--Mumford formulation of stability in terms of the weight of
limit points.

\subsubsection{Test configurations}\label{test-config}

Let $(M,\Omega)$ be a Hodge manifold, viewed as a polarized projective variety
with respect to a line bundle $L$ with $c_1(L)=\Omega/2\pi$.

\begin{defn} \cite{Do2}
A {\it test configuration} for $(M,L)$ is a polarized scheme $(X,\cE)$ over
$\C$ with a $\C\mult$ action $\alpha$ and a flat proper $\C\mult$-equivariant
morphism $p\colon X \to \C$ (where $\C\mult$ acts on $\C$ by scalar
multiplication) such that the fibre $(X_{t}= p^{-1}(t),\cE\restr{X_t})$ is
isomorphic to $(M,L)$ for some (hence all) $t\neq 0$.

$(X_0,\cE\restr{X_0})$ is called the {\it central fibre}.  Since $0 \in \C$ is
fixed by the action, $(X_0,\cE\restr{X_0})$ inherits a $\C\mult$ action, also
denoted by $\alpha$.

A test configuration is said to be a {\it product configuration} if $X=M\times
\C$ and $\alpha$ is given by a $\C\mult$ action on $M$ (and scalar
multiplication on $\C$).
\end{defn}
Since relevant properties of test configurations are unchanged if we replace
$\cE$ by $\cE^r$ for a positive integer $r$, we can let $\cE$ be a $\Q$-line
bundle in the definition above (i.e., $\cE$ denotes a `formal root' of a line
bundle $\cE^r$ for some positive integer $r$).

A particularly important class of test configurations are those associated to
a subscheme of $(M,L)$, as studied by Ross and Thomas~\cite{RT,RT2}. We shall
state it here for complex submanifolds of $(M,L)$, but the same definition
actually makes sense for subschemes.

\begin{defn}[\textit{Deformation to the normal cone}]
For a polarized complex manifold $(M,L)$, the {\it normal cone} of a complex
submanifold $Z$ is $\smash{\hat M}\cup_E P$, where $\smash{\hat M}$ is the
blow-up of $M$ along $Z$ with exceptional divisor $E=P(\nu_Z)$, $P=P(\cO\oplus
\nu_Z)$ and $\nu_Z$ is the normal bundle to $Z$ in $M$.  This is a singular
projective variety (for example the normal cone of a point $p\in \C P^1$ is
$\C P^1 \cup_p \C P^1$, which is a line-pair in $\C P^2$).

The normal cone is the central fibre of the family $p\colon X\to\C$ obtained
by blowing up $M\times\C$ along $Z\times\{0\}$ (where $p$ is the projection of
the blow-down to $\C$) called the {\it deformation to the normal cone} of $Z$
in $M$. We equip this with the polarization $\cE_{c} = \pi^*L
\otimes\cO(-cP)$, where $\cO(P)$ is the line bundle associated to the
exceptional divisor $P$, $\pi\colon X\to M$ is the projection of the blow-down
to $M$, and $c$ is a positive rational number such that $\cE_{c}$ is an ample
$\Q$-line bundle. This last condition gives an upper bound $\eps$ on $c$,
called the {\it Seshadri constant} of $Z$ with respect to $L$.

We let $\alpha$ be the $\C\mult$ action coming from the trivial action on $M$
and multiplication on $\C$. This clearly defines an action on $X$ with a lift
to $\cE_{c}$. Hence the deformation to the normal cone determines a family
of test configurations, parameterized by $c\in(0,\eps)\cap\Q$.
\end{defn}

\subsubsection{The Futaki invariant and K-stability}\label{futaki-stab}

K-stability is defined using a Hilbert--Mumford criterion, i.e., in terms of a
`weight' associated to each test configuration. This weight is given by the
Futaki invariant of the central fibre; however, since the latter is typically
a singular projective variety, we need an algebraic geometric definition of
the Futaki invariant. Such a definition has been given by Donaldson~\cite{Do2}.

Let $V$ be a scheme of dimension $n$ over $\C$ polarized by an ample line
bundle $L$ and suppose that $\alpha$ is a $\C\mult$ action on $V$ with a lift
to $L$.  Then $\alpha$ acts on the vector spaces $H_{k}=H^{0}(V, L^{k})$, $k
\in \Z^{+}$. If $w_{k}(\alpha)$ denotes the weight of the highest exterior
power of $H_{k}$ (that is, the trace $\mathrm{Tr}\,A_{k}$ of the infinitesimal
generator $A_{k}$ of the action) and $d_{k}$ denotes the dimension of $H_{k}$
then $w_{k}(\alpha)$ and $d_{k}$ are given by polynomials in $k$ for
sufficiently large $k$, of degrees at most $n+1$ and $n$ respectively.  For
sufficiently large $k$ the quotient $w_{k}(\alpha)/(kd_{k})$ can be expanded
into a power series with no positive powers. The {\it Futaki invariant}
$\mathfrak F(\alpha)$ is the residue at $k=0$ of this quotient, that is, the
coefficient of the $k^{-1}$ term in the resulting expansion. The Futaki
invariant is independent of the choice of lift of $\alpha$ to $L$.
\textup(When $V$ is a manifold, this definition coincides with Futaki's
original definition up to a normalization convention.\textup)

\begin{defn}
The Futaki invariant of a test configuration is defined to be the Futaki
invariant $\mathfrak F(\alpha)$ of the central fibre, where $\alpha$ denotes
the induced $\C\mult$ action.

A Hodge manifold $(M,L)$ is said to be K-{\it polystable} if the Futaki
invariant of any test configuration is nonpositive, and equal to zero if and
only if the test configuration is a product configuration.
\end{defn}
For the test configurations $(X,\cE_c)$ arising from a deformation to a normal
cone, Ross and Thomas~\cite{RT,RT2} show that the Futaki invariants $\mathfrak
F(\alpha_c)$ are rational in $c\in(0,\eps)\cap\Q$, where $\eps$ is the
Seshadri constant, and so can be extended to $c\in(0,\eps)$.

\begin{defn} A Hodge manifold $(M,L)$ is said to be {\it slope K-polystable}
if for the deformation to the normal cone of any nontrivial subscheme, the
Futaki invariant $\mathfrak F(\alpha_c)$ of the corresponding family
$(X,\cE_c)$ of test configurations is negative on $(0,\eps)$.
\end{defn}
Actually, the definition in~\cite{RT,RT2} is more subtle, since it requires
that $\mathfrak F(\alpha_\eps)<0$ unless $\eps$ is rational and the semi-ample
configuration $(X,\cE_\eps)$ is the pullback by a contraction of a product
configuration. We shall not need this refinement.

\subsection{Stable bundles and CSC K\"ahler metrics}\label{stab-CSC}

We now relate our results concerning CSC K\"ahler metrics on projective
bundles to stability theory for vector bundles. Recall that if $E\to S$ is a
holomorphic vector bundle over a compact k\"ahlerian $2d$-manifold
$(S,[\omega_h])$, the {\it slope} $\mu(E)$ is the number $\cb_1(E)\cdot
[\omega_h]^{d-1}$; $E$ is called (slope) {\it stable} or {\it semistable} if
$\mu(F) < \mu(E)$ or $\mu(F) \leq \mu(E)$ (respectively) for any proper
coherent subsheaf $F\subset E$; it is {\it polystable} if it is a direct sum
of stable vector bundles with the same slope; then, as is well-known, `stable'
$\Rightarrow$ `polystable' $\Rightarrow$ `semistable', and by the
Hitchin--Kobayashi correspondence, $E$ admits an Einstein--Hermitian
connection iff it is polystable.

There is a close analogy between K-stability for polarized K\"ahler manifolds
and slope stability for vector bundles. In particular, one might hope to find
a direct relation between the existence problem for CSC K\"ahler metrics on a
geometrically ruled complex manifold $P(E)$ over $S$ and the stability of
$E\to S$.  Notable progress in understanding the relation between
K-polystability of $P(E)$ and slope polystability of $E$ has been made by
Ross--Thomas~\cite{RT,RT2}, using their notion of slope K-polystability:
indeed if $F$ is a coherent subsheaf of $E$, then $P(F)$ is a subscheme of
$P(E)$ and deformation to the normal cone of $P(F)$ is a test configaration
which `destabilizes' $P(E)$ iff $\mu(F)>\mu(E)$ (see~\cite{RT}).

Using the general theory of CSC K\"ahler metrics, the work of~\cite{RT} shows
that if $E$ is not semistable with respect to an integral K\"ahler class
$[\omega_h]$ on $S$, then for all $k\gg 1$ the integral classes $2\pi
c_1(\cO(1)_E)+ k p^*[\omega_h]$ on $P(E)\stackrel{p}{\to} S$ do not contain
CSC metrics.  As a partial converse, Hong~\cite[Theorem A]{hong} shows that if
$E$ is polystable and $\mathfrak{h}_0(M)\to \mathfrak{h}_0(S)$ is surjective,
then there is a CSC metric in $2\pi c_1(\cO(1)_E) + k p^*[\omega_h]$, for each
$k\gg 1$, iff the Futaki invariant $\mathfrak F_\Omega$ vanishes.

To put our results in this context, let $P(E)\to S$ be admissible so
$E=E_0\oplus E_\infty$ for projectively-flat (and thus polystable) hermitian
vector bundles $E_0$ and $E_{\infty}$ with $\cb_1(E_{\infty}) - \cb_1(E_0) =
\sum_a[\omega_a/2\pi]$. Thus $E$ is determined up to tensor product by a line
bundle and is polystable iff it is semistable iff $\mu(E_0)=\mu(E_\infty)$.
With respect to a K\"ahler class $[\omega_h]=[\sum_a \omega_a/f_a]$ on $S$
(where $f_ax_a>0$), this condition reads
\begin{equation}\label{semistable}
0=\mu(E_{\infty}) - \mu(E_0) = \frac{(d-1)!}{2\pi} \Vol(S, {\textstyle\prod}_a
\omega_a/f_a)\Bigl({\textstyle\sum}_a d_af_a\Bigr),
\end{equation}
which can happen for suitable $f_a$ iff $\cb_1(E_{\infty}) - \cb_1(E_0)$ is
strictly indefinite; this is exactly the condition of Theorem~\ref{main-thm}
that ensures the existence of CSC metrics in a sufficiently small admissible
K\"ahler class $\Omega=\Xi+p^*[\sum_{a\in\cA}\omega_a/x_a]$, subject only to
the constraint that $\mathfrak F_\Omega(K)=0$. Note that $\Xi$ is equal to
$4\pi c_1(\cO(1)_E)$ up to a basic term (depending on the choice of
$E$)---this follows by integrating $\Xi^{ d_0+ d_{\infty}+1}$ over a fibre and
using the expression for $I(d_0,d_{\infty},0)$ from Appendix B.  Thus,
admissible classes play a similar role to those considered by Ross--Thomas and
Hong, and $k\gg 1$ corresponds to $|x_a|$ sufficiently small in our picture.
(We recall that this means that the fibres are small compared to the base.)
However, there is not a simple relation in general between those $\Omega$
containing a CSC metric and the $[\omega_h]$ with respect to which $E$ is
polystable: the approach of Ross--Thomas and Hong suggests taking
$f_a=x_a/(1+r_a x_a)$, for some $r_a$ depending only on $E$;
then~\eqref{semistable} agrees {\it asymptotically} with $\mathfrak
F_\Omega(K)=0$ in the limit $x_a\to 0$, but the two conditions define distinct
hypersurfaces in general.

Conversely, in the case $E_0$ and $E_\infty$ are line bundles over a product
$S=\Sigma_1\times\Sigma_2$ of two Riemann surfaces,
Theorem~\ref{no-CSC-product-surfaces} shows that polystability of $E$ with
respect to some K\"ahler class on $S$ (which is unique up to scale in this
case) is also necessary for the existence of a CSC metric in an admissible
K\"ahler class on $P(E)$.

Consider now the case that the base $S$ is a Riemann surface $\Sigma$ of genus
$\mathbf g$; the stability of a holomorphic vector bundle is then independent
of the choice of a K\"ahler class on $\Sigma$, and it is natural to
speculate~\cite{RT} that the notion of K-polystability of the projective
manifold $P(E)$ should be independent of the specific K\"ahler class, and to
conjecture that $P(E)$ admits a CSC K\"ahler metric iff $E$ is polystable. At
present (see~\cite{AT}) this conjecture is confirmed when $E$ is of rank $2$
(i.e., on geometrically ruled surfaces), when $\mathbf g\leq 1$ and $E$ is a
direct sum of line bundles (cf.~Corollary~\ref{rscase}---this always holds
when $\mathbf g=0$), or when $E$ is indecomposable and ${\mathbf g} \geq 2$.
Theorem~\ref{noCSCsigma} further confirms the conjecture in the case of
decomposable bundles of the form $E= E_0 \oplus E_{\infty}$ with $E_0$ and
$E_{\infty}$ polystable.

\subsection{Extremal K\"ahler metrics and relative K-polystability}

In recent work, G.~Sz\'ekelyhidi~\cite{szekelyhidi} has extended the theory of
K-polystability to cover extremal K\"ahler metrics, not just CSC K\"ahler
metrics. We briefly explain his ideas here.

\subsubsection{Motivation}\label{motiv} Recall that extremal
K\"ahler metrics are critical points for the $L_2$-norm of the scalar
curvature for metrics in a fixed K\"ahler class on a complex manifold
$(M,J)$~\cite{calabi1}. If we identify the K\"ahler class with a leaf of the
formal $G^c$ orbit described in \S\ref{infinitedim}, we are therefore looking
for critical points of $||\mu||^2$, where $\mu\colon X\to C^\infty_0(M,\R)$
and $X$ is the space of compatible complex structures on a compact symplectic
manifold $(M,\omega)$ with $H^1(M)=0$.

We can adapt the finite dimensional model of~\S\ref{finitedim} to this problem
by supposing that the Lie algebra $\mathfrak g$ is equipped with a
$G$-invariant inner product $\langle,\rangle$. Now, following
Sz\'ekelyhidi~\cite{szekelyhidi}, we note that the weight $w_x$ of the linear
action of the isotropy algebra $\mathfrak g_x$ on $\cL_x^*$ is given by
$w_x=\langle\beta_x,\cdot\rangle\colon \mathfrak g_x\to \C\mult$ for some
$\beta_x\in \mathfrak g_x$, which is the orthogonal projection of $\mu(x)$
onto $\mathfrak g_x$.  We refer to $\beta_x$ (or rather the induced vector
field on $X$) as the {\it extremal vector field}: for in the infinite
dimensional setting it agrees with the extremal vector field of Futaki and
Mabuchi (see \S\ref{futaki-invariant}).

Clearly $x$ is a critical point of $||\mu||^2$ iff $\beta_x$ is in $\mathfrak
g_x$. Using this, Sz\'ekelyhidi shows that $x$ is in the $G^c$ orbit of a
critical point of $||\mu||^2$ if and only if it is polystable {\it relative to
the extremal vector field}, i.e., for the action of the subgroup of $G^c$
whose Lie algebra is the subspace $\beta_x^\perp$ of the centralizer of
$\beta_x$. The Hilbert--Mumford criterion may then be modified as follows: the
{\it modified weight} $w_{x_0}(\alpha)-\langle \alpha,\beta_x\rangle
w_{x_0}(\beta_x)/\langle \beta_x,\beta_x\rangle$ of the limit point $x_0$
should be nonpositive for any one parameter subgroup $\alpha$ of the
centralizer of $\beta_x$, with equality if and only if $x_0=x$.

\subsubsection{The inner product and modified Futaki invariant}

Thus motivated, we return to the setting of \S\ref{futaki-stab} and define a
modified Futaki invariant of a polarized scheme $(V,L)$ (of dimension $n$ over
$\C$) relative to a $\C\mult$ action $\beta$. We first need to define an inner
product between such actions.

Assume then that $V$ has two $\C\mult$ actions $\alpha$ and $\beta$ with lifts
to $L$ and infinitesimal generators $A_k$ and $B_k$ of the actions on
$H_{k}$. Then for $k$ sufficiently large, $\mathrm{Tr}\,(A_k B_k)$ is a
polynomial of degree at most $n+2$. The inner product $\langle \alpha, \beta
\rangle$ is defined to be the coefficient of $k^{n+2}$ of the expansion of
$\mathrm{Tr}\,(A_k B_k) - {w_k(\alpha)w_k(\beta)}/d_k$ for large $k$, which is
independent of the lifts of $\alpha$ and $\beta$ to $L$: indeed it depends
only on the trace-free parts of $A_k$ and $B_k$. (When $V$ is a manifold, this
inner product coincides with Futaki--Mabuchi bilinear
form~\cite{futaki-mabuchi} up to a normalization convention.)

We define the {\it modified Futaki invariant} \cite{szekelyhidi} $\mathfrak
F_\beta(\alpha)$ of $\alpha$ relative to $\beta$ (assuming the action $\beta$
is nontrivial) to be
\begin{equation*}
\mathfrak F_\beta(\alpha) = \mathfrak F(\alpha) - \frac{\langle \alpha, \beta 
\rangle}{\langle \beta, \beta \rangle}\mathfrak F(\beta).
\end{equation*}

\subsubsection{Relative K-stability}

Let $(M,\Omega,L)$ be as in \S\ref{test-config} and suppose it has a
nontrivial $\C\mult$ action $\beta$.

\begin{defn} \cite{szekelyhidi}
A test configuration $(X,\cE)$ for $(M,L)$ is {\it compatible} with $\beta$ if
there is a $\C\mult$ action, also denoted by $\beta$, on $(X, \cE)$ preserving
$p\colon X\to\C$ and inducing the trivial action on $\C$, such that $\beta$
restricted to $(X_t,\cE\restr{X_t})$ coincides with the original action for $t
\neq 0$ under the isomorphism with $(M,L)$.

In this case we have an induced action on the central fibre $X_0$, also called
$\beta$, and the modified Futaki invariant of the test configuration is
defined to be $\mathfrak F_\beta(\alpha)$.

A polarized Hodge manifold $(M,L)$ with nontrivial $\C\mult$ action $\beta$ is
K-{\it polystable} {\it relative to} $\beta$ if the modified Futaki invariant
of any test configuration $(X,\cE)$ compatible with $\beta$ is nonpositive,
and equal to zero if and only if $(X,\cE)$ is a product.
\end{defn}

If $(M,L)$ has a $\C\mult$ action $\beta$ which preserves a subscheme $Z$, the
test configurations $(X,\cE_c)$ arising from the deformation to the normal
cone are compatible with $\beta$.  As in~\cite{RT,RT2}, $\mathfrak
F_\beta(\alpha_c)$ is rational in $c\in(0,\eps)\cap\Q$ and so extends to
$c\in(0,\eps)$. Thus, cf.~\S\ref{futaki-stab}, we have a notion of {\it slope}
K-polystability relative to $\beta$.

\begin{defn}
A polarized Hodge manifold $(M,L)$ with nontrivial $\C\mult$ action $\beta$ is
said to be {\it slope K-polystable} {\it relative to} $\beta$ if for the
deformation to the normal cone of any nontrivial subsheme preserved by
$\beta$, the modified Futaki invariant $\mathfrak F_\beta(\alpha_c)$ of the
corresponding family $(X,\cE_c)$ of test configurations is negative for
$c\in (0,\eps)$.
\end{defn}
As with the definition of (absolute) slope K-polystability, strictly speaking,
we should also require $\mathfrak F_\beta(\alpha_\eps)<0$ unless $\eps$ is
rational and $(X,\cE_\eps)$ is the pullback by a contraction of a product
configuration.

Nakagawa shows~\cite{nakagawa} that the \textup(Futaki--Mabuchi\textup)
extremal vector field associated to a Hodge K\"ahler manifold $(M,\Omega)$
with a maximal compact connected subgroup $G$ of $H_0(M)$ has closed orbits,
and therefore defines an effective $\C^{\times}$ action which we will refer to
as the {\it extremal} $\C^{\times}$ action of $(M, \Omega, G)$. Then, the
motivation of \S\ref{motiv} suggests the following
conjecture~\cite{szekelyhidi}.
\begin{conj}\label{stab-conj} Let $(M,\Omega, L)$ be a polarized Hodge
manifold and $G$ a maximal compact connected subgroup of $H_0(M)$. Then there
is a $G$-invariant extremal K\"ahler metric in $\Omega=2\pi c_1(L)$ if and
only if $(M,L)$ is K-polystable relative to the extremal $\C\mult$ action
of $(M,\Omega,G)$.
\end{conj}
As a motivating example, Sz\'ekelyhidi considers the deformation to the normal
cone of a the infinity section in a (higher genus) polarized geometrically
ruled surfaces $P(\cO \oplus \cL) \to\Sigma$. He finds that relative
K-polystability implies existence of extremal K\"ahler metrics of the type
constructed in \cite{christina1}, which are precisely the admissible metrics
on these bundles. In the next section we generalize this idea to arbitrary
admissible bundles. However, in doing so, we find that unless $\dim S\leq 4$,
we need to replace `K-polystable' by `slope K-polystable' in the above
conjecture.

\subsection{Relative K-polystability of admissible projective bundles}
\label{ext-exist}

We now consider the deformation to the normal cone $(X,\cE_{c},\alpha)$ of the
infinity section $e_{\infty} = z^{-1}(-1) = P(0 \oplus E_\infty)$ for an
admissible projective bundle $M=P(E_0\oplus E_\infty)\to S$ (with $\dim
S=2d$), polarized by a line bundle $L$ with $\Omega=2\pi c_{1}(L)$ admissible.

We therefore choose the admissible K\"ahler class $\Omega = \Xi + \sum_a
[\omega_a]/x_a$ to be integral (where $0<|x_a|\leq 1$ with equality iff
$a\in\{0,\infty\}$).  The Seshadri constant of this polarization is $2$, so we
take $c\in(0,2)\cap\Q$. Since the $\C\mult$ action $\beta$ induced by the
vector field $K$ preserves $Z$, $X$ is compatible.  We will use the letters
$\alpha,\beta$ to denote also the corresponding actions on the (polarized)
central fibre $(X_{0},L_0)$ and on the vector space $H^0(X_{0}, L_0^k)$, where
$L_0=\cE_{c}\restr{X_{0}}$.

Let us calculate the modified Futaki invariant of this configuration. For this
we first note that if $\cI_\infty \subset\cO_M$ is the ideal sheaf of
holomorphic functions vanishing on $e_\infty$, then for any $p\geq 0$,
$\cI^{p}_{\infty}/\cI^{p+1}_{\infty}$ is supported on $e_\infty$, and its
restriction is $S^p\nu_\infty^*$, where $\nu_\infty$ is the normal bundle to
$e_\infty$ in $M$.

Therefore, for $k$ sufficiently large, we have, as in~\cite{RT,szekelyhidi}
\begin{equation*}
H^{0}(X_{0}, L_0^{k}) = \bigoplus_{i=0}^{(2-c)k}
H^{0}(e_\infty,L|_{e_\infty}^k \otimes S^{2k-i}\nu_\infty^*)
\oplus \bigoplus_{j=1}^{ck}
H^{0}(e_\infty,L|_{e_\infty}^{k} \otimes S^{ck-j}\nu_\infty^*),
\end{equation*}
where $\alpha$ acts on the first direct sum with weight $0$ and on the
components of the second direct sum with weight $-j$. We can choose the lift
of $\beta$ to $L$ so that the weight of the induced action on
$H^0(e_\infty,L|_{e_\infty}^{k} \otimes S^{uk+v}\nu_\infty^*)$ is $(u-1)k+v$.

Now $S^p\nu_\infty^*$ is the direct image $q_* \cO(p)_{\nu_\infty}$, where
$\cO(-1)_{\nu_\infty}$ is the (fibrewise) tautological bundle of $q\colon \hat
e_\infty = P(\nu_\infty)\to e_\infty$. Also $\hat e_\infty$ may be identified
with $\smash{\hat S}$ via the obvious inclusion $i$ of $\smash{\hat S}$ into
$\smash{\hat M}=P(\cO\oplus\smash{\hat\cL})$ as the infinity section, and then
$i^*\cO(1)_{\nu_\infty}=\smash{\hat\cL}$. For convenience, we now drop the
hats, so that we have
\begin{align*}
H^{0}(X_{0}, L_0^{k})&=
\bigoplus_{i=0}^{(1-z)k} H^{0}(S,i^*L^k \otimes \cL^{2k-i})
\oplus \bigoplus_{j=1}^{(1+z)k} H^{0}(S,i^*L^k \otimes \cL^{(1+z)k-j})\\
&= \bigoplus_{i=0}^{2k} H^{0}(S,i^*L^k \otimes \cL^{2k-i}),
\end{align*}
where we have abused notation by writing $c-1=z$; {\it a priori} this has
nothing to do with the momentum map that we also denote by $z$, but notice
that it does take (rational) values in the same interval $(-1,1)$.  To compute
$d_{k}$, $\mathrm{Tr}\,A_{k}$, $\mathrm{Tr}\, B_{k}$, $\mathrm{Tr}\,
A_{k}B_{k}$, $\mathrm{Tr}\, B_{k}^{2}$, and thereby $\mathfrak
F_\beta(\alpha)$, we need only the dimensions of these vector spaces.  We note
that we only need to compute $d_{k}$, $\mathrm{Tr}\,A_{k}$ and $\mathrm{Tr}\,
B_{k}$ to subleading order in $k$, whereas for $\mathrm{Tr}\, A_{k}B_{k}$ and
$\mathrm{Tr}\, B_{k}^{2}$ the leading order term suffices.  Consequently we
will be dropping lower order terms without further comment. We also note that
since the Futaki invariant is defined in terms of ratios, we can ignore any
overall multiples.  Now by the Riemann--Roch formula and the ampleness of
$i^*L$ (in fact it is only semiample if $d_\infty>0$, but we can apply a
limiting argument in this case, as in~\cite{RT}), for sufficiently large $k$
we have that
\begin{align*}
h^{0}(S, i^*L^{k} \otimes \cL^{uk+v}) &= \chi(S, i^*L^{k} \otimes
\cL^{uk+v}) =\bigl(\mathrm{ch}(i^*L^{k} \otimes \cL^{uk+v})\cdot
\mathrm{Td}(S)\bigr)[S]\\ &= \bigl(c_{1}(i^*L^{k} \otimes
\cL^{uk+v}) + \tfrac 12 c_{1}(\cK^{-1}_{S})\bigr)^d[S]+O(k^{d-2})\\
&=
\biggl(\sum_a\frac{k+\bigl((u-1)k+v+s_a/2\bigr)x_a}{x_a}
[\omega_a/2\pi]\biggr)^d
[S]+O(k^{d-2})
\end{align*}
since the Ricci form of $S_a$ may be written $s_a\omega_a +\rho_{a,0}$ where
$\rho_{a,0}\wedge\omega_a^{d_a-1}=0$.  After an overall multiplication by
$(2\pi)^d/d!$ and $\prod_a x_{a}^{d_a}/\Vol(S_a,\omega_a)$, this is
\begin{equation*}
k^d p_{c}^{s}(u-1+{v}/{k})+O(k^{d-2}),
\end{equation*}
where $p_{c}^{s}(t) = \prod_a (1 + x_{a} (t+s_a/2k))^{d_{a}}$.  In order to
carry out the summations over $i$ and $j$ we use the trapezium rule, as
in~\cite[Lemma 4.7]{RT}.
\begin{lemma}
Let $f(x)$ be a polynomial and $b$ a rational number.  Then for
$\eps\in\{0,1\}$ and for $k\in \Z^+$ such that $bk$ is a positive integer, we
have
\begin{equation*}
\sum_{i=\eps}^{bk} f(i/k) = k \int_{0}^{b} f(t)\,dt + 
\frac{1}{2}(f(b) +(-1)^\eps f(0)) + O(k^{-1}).
\end{equation*}
\end{lemma}
The proof is easy (see e.g.~\cite{RT}): by linearity we can assume $f(x)=x^m$
and then use $\sum_{i=1}^N i^m = N^{m+1}/(m+1) + N^m/2 + O(N^{m-1})$ (which in
turn is an easy induction on $N$). We then obtain (up to an overall multiple),
that for any $r\geq 0$,
\begin{align*}
k^{-d-r} \mathrm{Tr}\,B_k^r & = k \int_{0}^{2k} (1-t)^r
p_{c}^{s}(1-t)\,dt+ \tfrac{1}{2}(p_{c}(1)+ (-1)^r p_{c}(-1))
+O(1/k) \\
& = k \alpha_{r} + \tfrac{1}{2} \beta_{r} + O(1/k),
\end{align*}
with $\alpha_r=\int_{-1}^1 p_c(t)t^r\,dt$ and $\beta_r$ as in~\eqref{betas}.
Setting $r=0$ gives $d_k$. Similarly, using the explicit
formula~\eqref{IVP-formula} for the extremal polynomial $F_\Omega(z)$ we
obtain
\begin{align*}
k^{-d-1}\mathrm{Tr}\,A_{k} & = k \int_{0}^{1+z} -t p_{c}^{s}(z - t)\,dt
- \tfrac{1}{2} (1+z) p_{c}(-1)+O(1/k)\\
& = -k \int_{-1}^{z} p_{c}(t)(z-t)\,dt - \frac{1}{2}\int_{-1}^z 
\biggl(\sum_{a} \frac{d_a s_a x_a}{1+x_a t}\biggr)
\Mpc(t) (z-t) dt\\
&\quad- \tfrac{1}{2} (1+z) p_{c}(-1)+O(1/k)\\
& =  - k \int_{-1}^{z} (z-t)p_{c}(t)\,dt\\
&\quad- \frac{1}{4} 
F_\Omega(z) + \frac{1}{4} \int_{-1}^{z}(At+B)p_{c}(t)(z-t)\,dt+O(1/k)
\displaybreak[0]\\
k^{-d-3}\mathrm{Tr}\,A_{k}B_{k} & = \int_{0}^{1+z}-t(z-t)p_{c}(z-t)\,dt
+O(1/k)\\
&= -\int_{-1}^{z}p_{c}(t)t(z-t)\,dt+O(1/k)
\end{align*}
where $A$ and $B$ are the solutions of~\eqref{system}.
Now we are ready to calculate $\langle \beta,\beta \rangle$, $\langle \alpha,
\beta \rangle$, $\mathfrak F(\beta)$, and $\mathfrak F(\alpha)$. (We omit the
dependence of $z$ for convenience.)
\begin{align*}
\langle \beta,\beta \rangle &= 
\frac{\alpha_{2}\alpha_{0} - \alpha_{1}^2}{\alpha_{0}}\\
\langle \alpha,\beta \rangle &= - \int_{-1}^{z} 
p_{c}(t)t(z-t)\,dt + \frac{\alpha_{1}}{\alpha_{0}}\int_{-1}^{z} 
p_{c}(t)(z-t)\,dt
\displaybreak[0]\\
\mathfrak F(\alpha) & = \mathrm{Res}_{k=0} \frac{(\mathrm{Tr}\,A_{k})_{1} + 
{(\mathrm{Tr}\,A_{k})_{0}}/{k}}{\alpha_{0}(1 + {\beta_{0}}/({2k\alpha_{0}}))}
= \frac{\alpha_{0}(\mathrm{Tr}\,A_{k})_{0} - 
\frac12{\beta_{0}(\mathrm{Tr}\,A_{k})_{1}}}{\alpha_{0}^{2}}\\
& = \frac{- \frac14{\alpha_{0}} 
F_\Omega(z) + \frac14{\alpha_{0}} \int_{-1}^{z}(At+B)p_{c}(t)(z-t)\,dt 
+\frac12{\beta_{0}}\int_{-1}^{z} p_{c}(t)(z-t)\,dt}{\alpha_{0}^{2}}\\
\mathfrak F(\beta) & = \mathrm{Res}_{k=0}\frac{\alpha_{1} + 
\beta_{1}/2k}{\alpha_{0}(1 + 
{\beta_{0}}/({2k\alpha_{0}}))}
= \frac{\beta_{1} \alpha_{0} - \beta_{0} \alpha_{1}}{2\alpha_{0}^{2}}
\end{align*}
where we have set $(\mathrm{Tr}\,A_{k})_{0}=- \frac{1}{4} F_\Omega(z) +
\frac{1}{4} \int_{-1}^{z}(At+B)p_{c}(t)(z-t)\,dt$ and
$(\mathrm{Tr}\,A_{k})_{1}= - k \int_{-1}^{z} p_{c}(t)(z-t)\,dt$.

Finally, we can calculate the modified Futaki invariant for our test 
configuration.
\begin{align*}
\alpha_{0}^{2}&\mathfrak F_\beta(\alpha)
=\alpha_{0}^{2}\bigl(\mathfrak F(\alpha) - 
{\langle \alpha, \beta \rangle \mathfrak F(\beta)}/{\langle \beta, \beta 
\rangle}\bigr) \\
& =  - \tfrac14{\alpha_{0}}
F_\Omega(z) + \tfrac14{\alpha_{0}} \int_{-1}^{z}(At+B)p_{c}(t)(z-t)\,dt 
+\tfrac12 {\beta_{0}}\int_{-1}^{z} p_{c}(t)(z-t)\,dt\\
&\quad +\frac{\alpha_{0}(\beta_{1}\alpha_{0} - 
\beta_{0}\alpha_{1})}{2(\alpha_{2}\alpha_{0} - \alpha_{1}^{2})}\biggl(
\int_{-1}^{z}  p_{c}(t)t(z-t)\,dt
-\frac{\alpha_{1}}{\alpha_{0}}\int_{-1}^{z}  p_{c}(t)(z-t)\,dt\biggr)
\displaybreak[0]\\
& = - \tfrac14{\alpha_{0}}
F_\Omega(z) + \tfrac14 {\alpha_{0}} \int_{-1}^{z}(At+B)p_{c}(t)(z-t)\,dt
- \tfrac14{\alpha_{0}}\int_{-1}^{z}At p_{c}(t)(z-t)\,dt\\
&\quad+ \frac{\beta_{0}(\alpha_{2}\alpha_{0} - \alpha_{1}^{2}) - 
\alpha_{1}(\beta_{1}\alpha_{0} - 
\beta_{0}\alpha_{1})}{2(\alpha_{2}\alpha_{0} - \alpha_{1}^{2})} 
\int_{-1}^{z} p_{c}(t)(z-t)\,dt
\displaybreak[0]\\
&= - \tfrac14{\alpha_{0}} 
F_\Omega(z) + \tfrac14{\alpha_{0}} \int_{-1}^{z}(At+B)p_{c}(t)(z-t)\,dt \\
&\quad - \tfrac14{\alpha_{0}}\int_{-1}^{z}At p_{c}(t)(z-t)\,dt
-\tfrac14{\alpha_{0}} \int_{-1}^{z} B p_{c}(t)(z-t)\,dt\\
&=-\tfrac14{\alpha_{0}} F_\Omega(z)
\end{align*}
which is a negative multiple of the extremal polynomial. It follows
immediately that if $(M,L)$ is slope K-polystable relative to $K=J\grad_g z$,
then $F_\Omega$ is positive on $(-1,1)$ and $\Omega$ contains an admissible
extremal metric by Theorem~\ref{main-thm}.

If $(M,L)$ is slope K-polystable in the absolute sense, then $\mathfrak
F(\alpha)$ is negative on $(-1,1)$ and hence nonpositive at $z=1$. Evaluating
the integrals in this case (and using $F_\Omega(1)=0$), we find that $A\geq
0$.  Now if we swap the roles of the zero and infinity sections (by
interchanging $E_0$ and $E_\infty$) then the analogous calculation shows that
$A\leq 0$ (we get the same formulae with the change of variables $z\mapsto
-z$). Thus $A=0$ and $\mathfrak F(\beta)=0$. (Intuitively, the reason we get
$\mathfrak F(\beta)=0$ is that in these limits, deformation to the normal cone
of the zero or infinity section is actually the pullback by a contraction of
the product configuration associated to $\pm\beta$, cf.~\cite{RT,RT2}.) Hence
$(M,L)$ is slope K-polystable relative to $\beta$, $F_\Omega$ is positive on
$(-1,1)$, and the admissible extremal metric is CSC.

This proves Theorem~\ref{stab=>ext}, providing evidence for the reverse
implication in Conjecture~\ref{stab-conj} (with relative K-stability replaced
by relative slope K-stability) because in our setting, the extremal vector
field is a nonzero multiple of $AK$.  This calculation also shows that the
forward implication in Conjecture~\ref{stab-conj} implies
Corollary~\ref{F-nonneg}, without referring to K-energy or the results of
Chen--Tian~\cite{CT,CT2}, providing further indirect evidence.  However, if we
use relative K-stability instead of relative slope K-stability, we can only
deduce that $F_\Omega$ is positive on $(-1,1)\cap\Q$ and hence nonnegative on
$(-1,1)$.  However, since $\Omega$ is integral, $F_\Omega$ has rational
coefficients, and so when $\dim S (=\sum_{a\in \cA} 2d_{a}) \leq 4$ it follows
that $F_{\Omega}$ is positive on $(-1,1)$: indeed
$F_{\Omega}(z)=(1+z)^{d_{0}+1}(1-z)^{d_{\infty}+1} Q(z)$ where $Q(z)$ is a
quadratic or cubic with rational coefficients, and the repeated roots of
such a polynomial must be rational.

On the other hand, the following examples show that positivity of the extremal
polynomial on $(-1,1)\cap \Q$ is not sufficient for the existence of an
extremal K\"ahler metric when $\dim S=6$.

\begin{ex} Let $S=\Sigma_1\times\Sigma_2\times\Sigma_3$ be a product of
hyperbolic Riemann surfaces $\Sigma_a$ with integral K\"ahler classes
$[\pm\omega_a]$. Then for any admissible projective line bundle $M$ over $S$
and any admissible integral K\"ahler class $\Omega$ on $M$ with parameters
$x_a\in\Q$, the extremal polynomial has the form
\begin{equation*}
F_\Omega(z) = (1-z^2)(p_c(z) + (1-z^2)(a_0+a_1z+a_2z^2)).
\end{equation*}
where the $a_j$ are determined by the constant gaussian curvatures $\pm s_a$
of $\Sigma_a$ (via~\eqref{extremal-deriv}--\eqref{extremal-values}).  However,
since we are free to choose the genera and degrees of the line bundles
defining $M$, the $s_a$ can be arbitrary rational numbers subject only to the
constraint that $s_a x_a<0$ (so that the gaussian curvatures are negative).
Hence we are free to choose the $a_j$ subject to this constraint.

We claim that for any rational $r>0$ and $x_1>x_2>0>x_3$, we can choose the
$a_j$ so that $F_\Omega(z)$ is a positive multiple of $(1-z^2)(z^2+r z-1)^2$
provided that
\begin{equation}\label{key-cond}
x_1 x_2 x_3 +x_1+x_2+x_3 =0.
\end{equation}
$F_\Omega$ then has a repeated root in $(0,1)$ and another in $(-\infty,-1)$
and for $r$ in an open subset of $\Q^+$, these roots are irrational. Obviously
for any $1>x_1>x_2>0$ rational,~\eqref{key-cond} has a unique rational
solution $x_3 = -(x_1+x_2)/(1+x_1 x_2)$ with $0>x_3>-1$, and it is elementary
to verify our claim by equating coefficients.  $F_\Omega''(z)$ is then
negative for large $z$ and has at least two roots in $(-1,1)$ and none in
$(1,\infty)$. We then check that for $r>8/5$, $F_\Omega''(-1)$ is negative and
so the other two roots (which must be real, since $F_\Omega$ must have four
inflection points) are in $(-\infty,-1)$.  Hence we can choose $x_1,x_2$ so
that $s_1 x_1<0$ and $s_2 x_2<0$, with $s_a$ defined
by~\eqref{extremal-deriv}--\eqref{extremal-values}. It automatically follows
that $s_3 x_3<0$.

These data then define a countably infinite family (parameterized by
$(x_1,x_2,r)$ in an open subset of $\Q^3$) of admissible projective line
bundles over products of three Riemann surfaces together with admissible
rational K\"ahler classes (which we can scale to be integral) such that
$F_\Omega$ is positive on $(-1,1)\cap\Q$, but has an irrational repeated
root in $(-1,1)$. By Theorem~\ref{main-thm} these K\"ahler classes then do not
contain an extremal K\"ahler metric.
\end{ex}

\appendix

\section{Relation to previous papers}\label{appA}

In this appendix we summarize the classification of compact K\"ahler
$2m$-manifolds $M$ with a hamiltonian $2$-form of order $\ell$ given
in~\cite[Theorem 2]{ACGT3}, and explain how Theorem~\ref{ACGTthm} follows from
this classification in the case $\ell=1$. We also give a nonexistence result
for extremal K\"ahler metrics when $\ell=2$.

\subsection{Summary of the classification}

Let $(M,g,J,\omega)$ be a compact connected K\"ahler $2m$-manifold with a
hamiltonian $2$-form $\phi$ of order $\ell$. Let $p(t)$ be the momentum
polynomial of $\phi$ and $K(t)=J\grad_g p(t)$ be the corresponding family of
hamiltonian Killing vector fields. We summarize results from~\cite{ACG2,ACGT3}
in italics.
\begin{fact}
The vector fields $\{K(t):t\in\R\}$ generate an effective isometric
hamiltonian action of an $\ell$-torus $\T$ on $M$ and $\Mp(t)$ has $m-\ell$
constant roots counted with multiplicity. This action is free on a connected
dense open subset $M^0$ of $M$.
\end{fact}
We let $S_\Delta$ be the stable quotient of $M$ by the induced action of the
complexified torus $\T^c$ and denote by $\eta_a$, for $a$ in a finite set with
$\leq m-\ell$ elements, the {\it distinct} constant roots of $\Mp(t)$ and by
$d_a$ their multiplicities.
\begin{fact}
$S_\Delta$ is covered by a product $\smash{\tilde S_\Delta}=\prod_a S_a$ of
K\"ahler $2d_a$-manifolds $(S_a,\pm g_a,\pm \omega_a)$, and $M^0\to S_\Delta$
is a principal $\T^c$-bundle.
\end{fact}
In~\cite{ACG2,ACGT3}, we took $a\in\{1,\ldots N\}$, but here we shall adopt
(in a moment) a different notation for the index set. We let
$\Mpc(t)=\prod_{a} (t-\eta_a)^{d_a}$ and write $\Mp(t)=\Mpc(t)\Mpn(t)$, where
$\Mpn(t)=\sum_{r=0}^\ell (-1)^r\sigma_r t^{\ell-r}$ and $\sigma_0=1$.  The
Killing vector fields $K_r:= J\grad_g \sigma_r$, for $r=1,\ldots \ell$, are
linearly independent on $M^0$.
\begin{fact}
The image $\Delta$ of the momentum map $(\sigma_1,\ldots\sigma_\ell)$ is a
simplex in $\mathfrak t^*\cong\R^{\ell}$, whose interior \textup(the image of
$M^0$\textup) is the image under the elementary symmetric functions of a
domain $D=\prod_{j=1}^{\smash\ell} (\beta_{j-1},\beta_j)$, where
$\beta_0<\beta_1<\cdots<\beta_\ell$. The roots of $\Mpn(t)$ define smooth,
functionally independent, pairwise distinct functions $\xi_j$
$(j=1,\ldots\ell)$ on $M^0$ which extend continuously to $M$ with image
$[\beta_{j-1},\beta_j]$. The codimension one faces of $\Delta$ may be labelled
$F_0,\ldots F_\ell$ such that on $F_j$, either $\xi_j=\beta_{j}$ or
$\xi_{j+1}=\beta_j$.
\end{fact}
The local description of the metric on $M^0$ is as follows.
\begin{fact}
There are $1$-forms $\theta_1,\ldots\theta_\ell$ on $M^0$ with
$\theta_r(K_s)=\delta_{rs}$ and $d\theta_r=\sum_a
(-1)^r\eta_a^{\ell-r}\omega_a$ and a function $\Theta$ of one variable
satisfying
\begin{gather}\label{Theta-pos}
(-1)^{\ell-j}\Theta>0\quad\text{on}\quad (\beta_{j-1},\beta_{j}),\\
\Theta(\beta_j) = 0,\qquad
\Theta'(\beta_j) = - \prod_{k\neq j}(\beta_{j}-\beta_k),
\label{Theta-bound}
\end{gather}
such that the K\"ahler structure on $M^0$ may be written
\begin{equation}\label{metric-ell}\begin{split}
g&=\sum_a \Mpn(\eta_a) g_a
+\sum_{j=1}^\ell \frac{\Delta_j}{\Theta(\xi_j)} d\xi_j^2
+\sum_{j=1}^\ell \frac{\Theta(\xi_j)}{\Delta_j}\Bigl(\sum_{r=1}^\ell
\sigma_{r-1}(\hat\xi_j)\theta_r\Bigr)^2,\\
\omega&=\sum_a \Mpn(\eta_a)\omega_a
+\sum_{r=1}^\ell d\sigma_r\wedge \theta_r,
\end{split}\end{equation}
where $\sum\nolimits_a \Mpn(\eta_a) g_a$ is the pullback of a local K\"ahler
product metric on $\smash{\hat S}$, $\Delta_j=\prod_{k \neq j} (\xi_j-\xi_k)$,
and $\sigma_r(\hat \xi_j)$ is the $r$th elementary symmetric function of
$\xi_1,\ldots \xi_\ell$ with $\xi_j$ omitted. \textup($\sigma_r$ itself is the
$r$th elementary symmetric function of $\xi_1,\ldots \xi_\ell$.\textup)
\end{fact}
The global description of $M$ in~\cite[Theorem 2]{ACGT3} was presented using
the blow-up $\smash{\hat M}$ of $M$ along the inverse image of the codimension
one faces $F_0,\ldots F_\ell$ of $S_\Delta$.
\begin{fact}
$\smash{\hat M}$ is $\T^c$-equivariantly biholomorphic to the $\C
P^\ell$-bundle $M^0\times_{\T^c}\C P^\ell\to S_\Delta$.
\end{fact}
The blow-up is encoded by fibrations $S_\Delta\to S_{F_j}$ for each $F_j$ (see
also~\cite[Proposition 6]{ACGT3}): either $S_\Delta=S_{F_j}$, or the fibration
is covered by the obvious projection $\smash{\tilde S_\Delta}\to\prod_{b\neq
a_j} S_b$ for some index $a_j$ such that $S_{a_j}$ is a complex projective
space and $(\pm g_{a_j},\pm \omega_{a_j})$ has constant holomorphic sectional
curvature $\pm \prod_{k\neq j}(\beta_{j}-\beta_k)$

We unify these cases here by introducing, if $S_\Delta=S_{F_j}$, an additional
index $a_j$ with $d_{a_j}=0$ and $S_{a_j}=\C P^0$ (a point).  We denote the
new index set by $\smash{\hat\cA}$ and take $a\in\smash{\hat\cA}$ unless
otherwise stated: the additional indices make no difference to the previous
formulae.  We still have $\smash{\tilde S_\Delta}=\prod_a S_a$, and now for all
$F_j$, $S_\Delta\to S_{\smash{F_j}}$ is a $\C P^{d_{a_j}}$-bundle covered by
$\smash{\tilde S_\Delta}\to\prod_{b\neq a_j} S_b$. The map $j\to a_j$ is
injective~\cite{ACGT3} and so $\smash{\hat\cA}$ is the union of a set $\cA$
and the injective image of $\{0,\ldots\ell\}$ (under $j\mapsto a_j$).
\begin{fact}
For $a\in\cA$, either $\eta_a<\beta_0$ or $\eta_a>\beta_\ell$, according to
the sign of $(\pm g_a,\pm\omega_a)$, whereas for $j=\{0,\ldots\ell\}$,
$\eta_{a_j}=\beta_j$.
\end{fact}
The formula~\eqref{metric-ell} for the metric on $M^0$ leads to a
description~\cite{ACGT3} of $\smash{\hat M}$ as a projective bundle
$P(\cL_0\oplus \cL_1\oplus\cdots \oplus \cL_\ell)\to S_\Delta$ together with
formulae for the Chern classes of $\cL_j$ on the covering $\smash{\tilde
S_\Delta}$. To obtain instead a description of $M$, we need one further
ingredient, which follows easily by considering the form of the covering
transformations and the fact that $\smash{\tilde S_\Delta}\to\prod_{b\neq a_j}
S_b$ covers the fibration $S_\Delta\to S_{\smash{F_j}}$.
\begin{lemma} The projection $\smash{\tilde S_\Delta} \to \smash{\tilde S}
:=\prod_{a\in\cA} S_a$ descends to realize $S_\Delta$ as a fibre product of
flat projective unitary $\C P^{d_{a_j}}$-bundles over a quotient $S$ of
$\tilde S$.
\end{lemma}
An important class of flat projective unitary $\C P^r$-bundles on $S$ are
those of the form $P(E)$, where $E$ is a rank $r+1$ projectively-flat
hermitian vector bundle on $S$. If $S$ is simply connected, then any flat
projective unitary $\C P^r$-bundle is trivial, hence of the form $P(E)$ with
$E\cong{\mathcal E}\otimes \C^{r+1}$ for a holomorphic line bundle ${\mathcal
E}$.  In general the obstruction to the existence of $E$ is given by a torsion
element of $H^2(S, \cO^*)$ (cf.~\cite{elencwajg-narasimhan}). In particular,
such an $E$ always exists if $S$ is a Riemann surface.

Let us suppose that $S_\Delta= P(\vE_0)\times_S P(\vE_1)\times_S\cdots\times_S
P(\vE_\ell)\to S$, where each $\vE_j\to S$ is projectively-flat hermitian of
rank $d_j+1$.  We are free to choose the $\vE_j$ so that $\smash{\hat
M}=P\bigl(\cO(-1)_{E_0}\oplus \cO(-1)_{E_1}\oplus\cdots
\oplus\cO(-1)_{E_\ell}\bigr)$ where $\cO(-1)_{E_j}$ is the (fibrewise)
tautological line bundle over $P(\vE_j)$ (trivial over the other factors of
$S_\Delta$).  From the description of the blow-up in~\cite{ACGT3}, we
immediately deduce the following (in which we write $\cb_1(E)=c_1(E)/\rk E$).
\begin{fact}
$M$ is $\T^c$-equivariantly biholomorphic to
$P(\vE_0\oplus\vE_1\cdots\oplus\vE_\ell)\to S$ and for any $i\neq j$,
$\cb_1(E_j)-\cb_1(E_i)=\tfrac 1 2 \sum_a\bigl(\prod_{k\neq i} (\eta_a-\beta_k)
- \prod_{k\neq j} (\eta_a-\beta_k) \bigr) [\omega_a/2\pi]$.
\end{fact}

\subsubsection*{Derivation of Theorem~\ref{ACGTthm}}

In order to derive Theorem~\ref{ACGTthm} from the above, it suffices to
rescale $g$ so that we can take $\beta_0=-1$ and $\beta_1=1$. Then we set
$\eta_a=-1/x_a$ and change the sign of $\omega_a$ for all $a$. We also write
$\smash{\hat S}$ for $S_\Delta$, and replace the index set $\{0,1\}$ by
$\{0,\infty\}$ so that we can take $\smash{\hat\cA}= \{0,\infty\}\cup\cA$
where $\cA$ is a finite subset of $\Z^+$, but these changes are purely
cosmetic.

\subsection{A nonexistence result for order 2 extremal K\"ahler metrics}

In this paper we study only hamiltonian $2$-forms of order $1$. As a partial
justification for this restriction, we now consider the lowest interesting
dimension for the order $2$ case, and show that any extremal metric on a
compact K\"ahler $6$-manifold compatible with a hamiltonian $2$-form of order
$2$ is a Fubini--Study metric on $\C P^3$.

In this situation, the momentum polynomial has nonconstant roots $\xi_{1}$ and
$\xi_{2}$ and one constant root $\eta$ so $\#\cA\leq 1$ and $\Mp(t) =
(t-\eta)(t-\xi_{1})(t-\xi_{2})$. The stable quotient $\Sigma$ of $(M,J)$ by
the complexified $\T^c$ action is a compact Riemann surface with K\"ahler
structure $(g_{\Sigma},\omega_\Sigma)$.

We can set $\beta_0=-1$ and $\beta_2=1$ and write $\beta_1=\beta$ (where
$|\beta|<1$). If $\cA$ is empty, $(M,J)$ is biholomorphic to $\C P^{3}$;
otherwise $|\eta|>1$ and $(M,J)$ is $\T^c$-equivariantly biholomorphic to
$M=P(\cL_0 \oplus \cL_1 \oplus \cL_2) \to \Sigma$, where $\cL_j$ are
holomorphic line bundles on $\Sigma$ such that (without loss) $\cL_1$ is
trivial and
\begin{equation}\label{planebdlsintegrality} \begin{split}
c_1(\cL_0)&=\tfrac 12 (\eta -1)(\beta+1)[\omega_\Sigma/2\pi],\\
c_1(\cL_2)&=\tfrac 12 (\eta +1)(\beta-1)[\omega_\Sigma/2\pi].
\end{split}\end{equation}
The K\"ahler metric on $M$ is determined by a function $\Theta(t)$ satisfying
positivity and boundary conditions which imply that $\Theta(t)=F(t)/(t-\eta)$
where $F(t)=H(t)((t-\eta)+H(t)Q(t))$ for some function $Q(t)$, and
$H(t)=(1-t^{2})(t-\beta)$.

If $g$ is extremal and the extremal vector field is tangent to the fibres of
$M\to\Sigma$, then by~\cite{ACG2}, $F(t)$ is a polynomial of degree at most
$5$ and $g_{\Sigma}$ has scalar curvature $-F''(\eta)$. This forces $Q(t)=0$
and so the scalar curvature of $g_{\Sigma}$ is $2( 3 \eta^2-2\beta\eta-1)$
which is positive since $|\eta | > 1$ and $|\beta| < 1$. Hence $\Sigma=\C
P^1$.  Since $\frac14 {\Scal_{g_\Sigma}}[\omega_\Sigma/2\pi]$ is a primitive
integral class, \eqref{planebdlsintegrality} implies that
\begin{equation*}
(\eta\mp1)(\beta\pm1) = q^{\pm} ( 3 \eta^2-2\beta\eta-1 )
\end{equation*}
for some nonzero integers $q^{\pm}$. We remark that these formulae show that
the relation between $q^\pm$ and $(\eta,\beta)$ is birational, in fact the
restriction to $\R^2$ of a quadratic transformation of $\C P^2$. In any case,
$\eta$ is constant on the lines through $(q^+,q^{-})=(1,1)$, and $\beta=\pm1$
on the lines $q^\pm=0$ and $2q^{\pm} - q^{\mp}=1$, the latter being the lines
on which $\eta=\pm1$. It follows straightforwardly that $|\eta|>1$ and
$|\beta|<1$ iff $q^+ > 0$, $q^- < 0$ and $2 q^+ - q^{-} < 1$ or vice-versa
(swap plus and minus)---which is impossible as $|q_\pm|\geq 1$.  We therefore
have the following nonexistence result.
\begin{thm}
\label{nonexistenceplanebundles}
A compact extremal K\"ahler $6$-manifold $(M,J, g,\omega)$ which admits a
hamiltonian $2$-form of order $2$ with the extremal vector field tangent to
the $\T^c$-orbits is isometric to $\C P^{3}$ with a Fubini--Study metric.
\end{thm}

\section{Proof of Lemma \ref{asymptotics}} \label{appB}

In this appendix we prove Lemma~\ref{asymptotics} by computing the asymptotics
as $x_a\to 0$, for $a\in\cA$, of the solution $(A,B)$ of the
system~\eqref{system}, i.e., $A \alpha_1 + B \alpha_0 =-2\beta_0$, $A \alpha_2
+ B \alpha_1 =-2\beta_1$, where $\alpha_r=\int_{-1}^1 \Mpc(t) t^r dt$ and
$\beta_r$ are as in~\eqref{betas}.  In order to do this, we rewrite $\beta_0$
and $\beta_1$ as integrals using the obvious identities
\begin{align*}
\Mpc(1)+\Mpc(-1) &= \int_{-1}^1 \frac{d}{dt} \bigl(t \Mpc(t)\bigr) dt
= \int_{-1}^1 \biggl(1+\sum_a \frac{d_a x_at}{1+x_at}\biggr) \Mpc(t)dt\\
\Mpc(1)-\Mpc(-1) &= \int_{-1}^1 \frac{d}{dt} \bigl(t^2 \Mpc(t)\bigr) dt
= \int_{-1}^1 \biggl(2+\sum_a \frac{d_a x_at}{1+x_at}\biggr) \Mpc(t)t\,dt
\end{align*}
to obtain
\begin{align*}
\beta_0&= \int_{-1}^1 (1+t)^{d_0}(1-t)^{d_\infty}\biggl( 1+d_0+d_\infty
+\frac{d_0^2}{1+t}+\frac{d_\infty^2}{1-t}+\sum_{a\in\cA}
\frac{d_a x_a(s_a+t)}{1+x_a t}\biggr)\\&\qquad\times
\biggl(\;\prod_{a\in\cA}(1+x_at)^{d_a}\biggr) dt
\displaybreak[0]\\
\beta_1&= \int_{-1}^1 (1+t)^{d_0}(1-t)^{d_\infty}\biggl( 2+d_0+d_\infty
+\frac{d_0^2}{1+t}+\frac{d_\infty^2}{1-t}+\sum_{a\in\cA}
\frac{d_a x_a(s_a+t)}{1+x_a t}\biggr)\\&\qquad\times
\biggl(\;\prod_{a\in\cA}(1+x_at)^{d_a}\biggr) t\,dt.
\end{align*}
The asymptotics of $\alpha_0,\alpha_1,\alpha_2,\beta_0$ and $\beta_1$ are
given by integrals of the form
\begin{equation*}
I(m,n,k) = \int_{-1}^1 (1+t)^{m}(1-t)^{n} t^k\, dt.
\end{equation*}
Integrating by parts and using $2I(m,n,k+1) = I(m+1,n,k)-I(m,n+1,k)$,
\begin{align*}
I(m,n,0)&= \frac{2^{m+n+1}\,m!\, n!}{(m+n+1)!},\qquad
I(m,n,1)= \frac{2^{m+n+1} (m-n)\, m!\, n!}{(m+n+2)!},\displaybreak[0]\\
I(m,n,2)&= \frac{2^{m+n+1} (m^2+n^2+m+n-2mn+2)\,m!\, n!}{(m+n+3)!}.
\end{align*}
These are rather complicated, so we manipulate the integrals using the
identities
\begin{align*}
I(m-1,n,0)m^2+I(m,n-1,0)n^2&=\tfrac12 I(m,n,0)(m+n+1)(m+n)\\
I(m-1,n,1)m^2+I(m,n-1,1)n^2&=\tfrac12 I(m,n,1)(m+n-1)(m+n+2)\\
I(m-1,n,2)m^2+I(m,n-1,2)n^2&=\tfrac12 I(m,n,2)(m+n+3)(m+n)\\
&\qquad-I(m,n,1)(m-n)
\end{align*}
and thus obtain, up to $O(x^2)$,
\begin{align*}
\alpha_k &= I(d_0,d_\infty,k)+I(d_0,d_\infty,k+1)
\sum_{a\in\cA} d_ax_a,\\
\beta_0 &= \tfrac12 I(d_0,d_\infty,0)(1+d_0+d_\infty)(2+d_0+d_\infty)
+I(d_0,d_\infty,0)\sum_{a\in\cA} d_as_ax_a\\
&+\tfrac12 I(d_0,d_\infty,1)(1+d_0+d_\infty)(2+d_0+d_\infty)
\sum_{a\in\cA} d_ax_a,\displaybreak[0]\\
\beta_1 &= \tfrac12 I(d_0,d_\infty,1)(1+d_0+d_\infty)(2+d_0+d_\infty)
+I(d_0,d_\infty,1)\sum_{a\in\cA} d_as_ax_a\\
&+\bigl(\tfrac12 I(d_0,d_\infty,2)(3+d_0+d_\infty)(2+d_0+d_\infty)
-I(d_0,d_\infty,1)(d_0-d_\infty)\bigr)\sum_{a\in\cA} d_ax_a.
\end{align*}
Direct computation with these formulae and the identity
$I(m,n,1)(m+n+2)=I(m,n,0)(m-n)$ now shows that, up to $O(x^2)$,
\begin{align*}
\frac{\alpha_0\beta_1-\alpha_1\beta_0}{\alpha_0\alpha_2-\alpha_1^2}
&= (2+d_0+d_\infty)\sum_{a\in\cA} d_ax_a\\
\frac{\alpha_2\beta_0-\alpha_1\beta_1}{\alpha_0\alpha_2-\alpha_1^2}
&= \tfrac12(1+d_0+d_\infty)(2+d_0+d_\infty)
+ \sum_{a\in\cA} d_as_ax_a+ (d_\infty-d_0) \sum_{a\in\cA} d_ax_a.
\end{align*}
Multiplying by $-2$ completes the proof.

\end{document}